\documentclass[10pt]{article}
\usepackage{amsfonts,latexsym,amsmath,amsthm,amssymb}
\usepackage{mathrsfs}

\newcommand{\D}{\Delta}
\newcommand{\p}{\partial}
\renewcommand{\phi}{\varphi}

\begin{document}

\newtheorem{theorem}{\indent Theorem}[section]
\newtheorem{proposition}[theorem]{\indent Proposition}
\newtheorem{definition}[theorem]{\indent Definition}
\newtheorem{lemma}[theorem]{\indent Lemma}
\newtheorem{remark}[theorem]{\indent Remark}
\newtheorem{corollary}[theorem]{\indent Corollary}

\begin{center}
{\large \bf Non-autonomous dynamics of wave equations \\
with nonlinear damping and critical nonlinearity} \vspace{0.5cm}\\
Chunyou Sun$^{1,2}$, Daomin Cao$^{1}$\\
{\small \it $^1$Institute of Applied Mathematics, Chinese Academy of Sciences\\
Beijing, 100080, China} \\
{\small \it $^2$School of Mathematics and statistics, Lanzhou
University\\ Lanzhou, 73000, China}\\
{\small E-mails: cysun@amss.ac.cn;  dmcao@amt.ac.cn}\\ \vspace{0.3 cm} and \\
Jinqiao Duan\\
{\small \it Department of Applied Mathematics, Illinois Institute
of Technology \\Chicago, IL 60616, USA \\
E-mail:  duan@iit.edu}
\end{center}

\numberwithin{equation}{section}
\footnote[0]{\hspace*{-7.4mm}
Date: March 27, 2006  submitted; September 7, 2006 revised. 
A part of this work was done while J. Duan was visiting the American
Institute of Mathematics, Palo Alto, California, USA. This work was
partly supported by the NSF Grants DMS-0209326 \& DMS-0542450, the
NSFC Grant 10601021 and the Outstanding Overseas Chinese Scholars
Fund of the Chinese Academy of Sciences. }


\begin{abstract}
The authors consider non-autonomous dynamical behavior of
wave-type evolutionary equations with nonlinear damping and
critical nonlinearity. These type of waves equations are
formulated as   non-autonomous dynamical systems (namely,
cocycles). A sufficient and necessary condition for the existence
of pullback attractors is established for norm-to-weak continuous
non-autonomous dynamical systems, in terms of pullback asymptotic
compactness or pullback $\kappa-$contraction criteria. A technical
method for verifying pullback asymptotic compactness, via
contractive functions, is devised. These results are then applied
to the wave-type evolutionary equations with nonlinear damping and
critical nonlinearity, to obtain the existence of pullback
attractors. The required pullback asymptotic compactness for the
existence of pullback attractors is fulfilled by some new a priori
estimates for concrete wave type equations arising from
applications. Moreover, the pullback $\kappa-$contraction
criterion for the existence of pullback attractors is of
independent interest.


\centerline{\emph{Dedicated to Philip Holmes on the occasion of his
60th birthday}}
\end{abstract}

\section{Introduction}

Nonlinear wave phenomena occur in various systems   in physics,
engineering, biology and geosciences \cite{Ba,CV,Ha,Te,Majda,
Pedlosky}. At the macroscopic level, wave phenomena   may be
modeled by hyperbolic wave type partial differential equations. We
consider the following non-autonomous wave equations with
nonlinear damping, on a bounded domain $\Omega$ in $\mathbb{R}^3$,
with smooth boundary $\partial \Omega$:
\begin{equation}\label{1.1}
u_{tt}+h(u_t)-\Delta u+f(u,t)=g(x,\,t)\quad x\in\Omega
\end{equation}
subject to the boundary condition
\begin{equation}\label{1.2}
u |_{\partial \Omega}=0,
\end{equation}
and the initial conditions
\begin{equation}\label{1.3}
u(x,0)=u_{0}(x),\quad u_t(x,0)=v_0(x).
\end{equation}
Here $h$ is the nonlinear damping function, $ f$ is the
nonlinearity, $ g$ is a given external time-dependent forcing, and
$\D = \p_{x_1x_1}+\p_{x_2x_2} + \p_{x_3x_3}   $ is the Laplace
operator.

Equation \eqref{1.1} arises as an evolutionary mathematical model in
various systems. For example, (i) modeling a continuous Josephson
junction with specific $h, g$ and $ f $ \cite{LSC}; (ii) modeling a
hybrid system of nonlinear waves and nerve conduct; and (iii) when
$h(u_t)=ku_t$ and $f(u)=|u|^ru$, the equation \eqref{1.1}  models a
phenomenon in quantum mechanics \cite{Cha,CY,GM,Te}.

For the autonomous case of \eqref{1.1}, i.e., when $f$ and $g$ do
not depend on time $t$ explicitly, the asymptotic behaviors of the
solutions have been studied extensively in the framework of global
attractors; see, for example, \cite{ACH,BV,Ba,CV,Ha} for the
linear damping case, and \cite{CL1,CL2,CL3,Fe,SYZ1} for the
nonlinear damping case.

In this paper, we consider the non-autonomous case, especially
when the damping $h$ is   nonlinear  and when the nonlinearity $f$
has critical exponent (see below). For a non-autonomous dynamical
system  like \eqref{1.1}-\eqref{1.3}, the solution map does not
define a semigroup and instead, it defines a two-parameter
\emph{process}, or \emph{cocycle}. Pullback attractors are
appropriate geometric objects for describing asymptotic dynamics
for cocycles. We will briefly introduce basic concepts for
non-autonomous dynamical systems in \S 3.
We will discuss the asymptotic dynamics of \eqref{1.1}-\eqref{1.3}
via   pullback attractors of the corresponding cocycle. This
dynamical framework allows us to handle more general
non-autonomous time-dependency; for example, the external force
$g$ needs to be neither almost periodic nor translation compact in
time.

Our basic assumptions about nonlinear damping $h$, nonlinearity
$f$ and forcing $g$ are as follows. Let $g(x,t)$ be in the space $
L_{loc}^2(\mathbb{R};L^2(\Omega))$, of locally square-integrable
functions, and assume that the   functions $h$ and $f$ satisfy the
following conditions:
\begin{equation}\label{1.4}
h\in \mathcal{C}^1(\mathbb{R}),\quad h(0)=0,\quad h~\text{\rm
strictly increasing},
\end{equation}
\begin{equation}\label{1.5}
\liminf_{|s|\to \infty}h'(s)>0,
\end{equation}
\begin{equation}\label{1.6}
|h(s)|\leq C_1(1+|s|^p),
\end{equation}
where $p \in [1,\, 5)$ which will be given precisely later; $f\in
\mathcal{C}^1(\mathbb{R}\times \mathbb{R};\,\mathbb{R})$ and
satisfies
\begin{equation}\label{1.7}
F_s(v,s)\leqslant \delta^2 F(v,s) +C_{\delta},\quad F(v,s)\geqslant
-mv^2-C_m,
\end{equation}
\begin{equation}\label{1.8}
|f_v(v,s)|\leq C_2(1+|v|^{q}),~|f_s(v,s)|\leq C_3(1+|v|^{q+1}),
\end{equation}
\begin{equation}\label{1.9}
f(v,s)v-C_4F(v,s)+mv^2\geqslant -C_m,\quad \forall ~(v,s)\in
\mathbb{R}\times \mathbb{R},
\end{equation}
where $0\leqslant q\leqslant 2$, $F(v,s)=\int_0^vf(w,s)dw$ and
$\delta$, $m$ are sufficiently small which will be determined in
$Lemma$ \ref{l5.3}. The number $q=2$ is called the \emph{critical
exponent}, since the nonlinearity $f$ is not compact in this case
(i.e., for a bounded subset $B\subset H_0^1(\Omega)$, in general,
$f(B)$ is not precompact in $L^2(\Omega)$). This is an essential
difficulty in studying the asymptotic behavior even for the
autonomous cases \cite{ACH,BV,Ba,CL1,CL2,CL3,Fe,SYZ1}. The
assumptions \eqref{1.4}-{\eqref{1.6} on $h$ are similar to those in
\cite{CL3,Fe,Kh,SYZ1} for the autonomous cases, while the assumption
$1 \leq p<5$ is due to the need for estimating
$\int_{\Omega}g(u_t)u$ by $\int_{\Omega}g(u_t)u_t$ and
$\int_{\Omega}|\nabla u|^2$ via Sobolev embedding. Finally, the
assumptions \eqref{1.7}-\eqref{1.9} are similar to the conditions
used in Chepyzhov \& Vishik \cite{CV} for non-autonomous cases
 but with linear damping.

\medskip

Let us recall some recent relevant research in this area.

The existence of pullback attractors are established for the
strongly dissipative non-autonomous dynamical systems such as those
generated by parabolic type partial differential equations, e.g.,
the non-autonomous 2D Navier-Stokes equation and some non-autonomous
reaction diffusion equations; see
\cite{ARS,CLR,CLR1,CR,Ch,CD,CKS,LLR} and the references therein.
However, the situation for the hyperbolic wave type systems  is less
clear. For the linear damping case $h(v)=kv$ with  a constant $k>0$
and $q<2$ (subcritical), Chepyzhov \& Vishik \cite{CV} have obtained
the existence of a uniform absorbing set when $g$ is translation
bounded in time (i.e., $g\in L^2_b(\mathbb{R};L^2(\Omega))$), and
the existence of a uniform attractor when $g$ is translation compact
in time (i.e., $g\in L^2_c(\mathbb{R};L^2(\Omega))$).

Under the assumptions that $g$ and $\p_t g$ are both in  the space
of bounded continuous functions
$\mathcal{C}_b(\mathbb{R},L^2(\Omega ))$, $h$ has bounded positive
derivative, and furthermore, $f$ is of critical growth (i.e.,
$q=2$), Zhou \& Wang \cite{ZW} have proved the existence of kernel
sections and obtained uniform bounds of the Hausdorff dimension of
the kernel sections. Caraballo et al. \cite{CKR} have discussed
the pullback attractors for the cases of linear damping and
subcritical nonlinearity ($q<2$) .


 As in the autonomous case, some kind of compactness of the
cocycle is a key ingredient for the existence of pullback
attractors of cocycles. The corresponding compactness assumption
in Cheban \cite{Ch} is that the cocycle   has a compact attracting
set. Recently, Caraballo et al \cite{CLR} have established a
criterion for the existence of pullback attractors via pullback
asymptotic compactness, and illustrated their results with the 2D
Navier-Stokes equation.

For the autonomous \emph{linearly} damped wave equations, Ball
\cite{Ba} proposed a     method to verify the   asymptotic
compactness for the corresponding solution semigroup. This
so-called energy method has been generalized \cite{LS, MRW}
  to some non-autonomous cases. However, for our problem, due
to the nonlinear damping, it appears   difficult to apply the
method  of Ball \cite{Ba}. Moreover, a decomposition technique
\cite{ACH,CV,Fe,Ha,PZ,SY} has been successfully applied   to
verify the asymptotic smoothness  of the corresponding solution
semigroup for autonomous wave equations.

In this paper, after some preliminaries, we first introduce the
\emph{pullback $\kappa-$contraction} concept, a generalization of
$\kappa-$contraction from autonomous systems to non-autonomous
systems.  Then we establish a  criterion for the existence of
pullback attractors, in terms of pullback $\kappa-$contraction or
pullback asymptotic compactness. This criterion is for a class of
``weakly" continuous cocycles (i.e., the so-called norm-to-weak
continuous cocycles; see \S 3 below). Thirdly, we show that the
pullback $\kappa-$contraction is not equivalent to the pullback
asymptotic compactness,  unless the cocycle mapping has a nested
bounded pullback absorbing set (see $Definition$ \ref{d3.2}
below). This fact   is different from the autonomous semigroup
cases. Moreover, we propose a   technique  for verifying pullback
asymptotic compactness. Finally, we apply these results to show
the existence of pullback attractors for the non-autonomous
hyperbolic wave system   \eqref{1.1}-\eqref{1.3}.


Due to the difference between the cases $p=1$ and $1<p<5$ for the
nonlinear damping exponent $p$, we propose the following two kinds
of assumptions.
\begin{enumerate}
\item[] \hspace{-0.5 cm}{\bf Assumption  I.}

$h$ satisfies \eqref{1.4}-\eqref{1.6} with $p=1$ and there is a
$C_0$ such that $C_0|u-v|^2\leqslant (h(u)-h(v))(u-v)$;

$g$ satisfies
\begin{equation}\label{1.10}
\int_{-\infty}^{t}e^{\beta s}\int_{\Omega}|g(x,s)|^2dxds < \infty
~~\text{for each}~~t\in \mathbb{R},
\end{equation}
where $\beta(<C_0)$ is constant depending on the coefficients of $h$
and $f$, which will be determined in the proof of $Lemma$
\ref{l5.3};

$f$ satisfies \eqref{1.7}-\eqref{1.9}.

\item[] \hspace{-0.5 cm}{\bf Assumption  II.}

$h$ satisfies \eqref{1.4}-\eqref{1.6} with $1 \leqslant p <5$; And
\begin{equation}\label{1.11}
g \in L^{\infty}(\mathbb{R},\,L^2(\Omega));
\end{equation}
In addition to \eqref{1.7}-\eqref{1.9}, $f$ satisfies also
\begin{equation}\label{1.12}
F_s(v,s) \leqslant 0 \quad \text{for all}~~ (v,s)\in
\mathbb{R}\times \mathbb{R}.
\end{equation}
\end{enumerate}

We remark that the technical hypotheses \eqref{1.11} and
\eqref{1.12} in {\it Assumptions $II$} are mainly for the existence
of pullback absorbing set; see $Lemma$ \ref{l5.3} below or Haraux
\cite{Ha2} for more details. Our method  for verifying the
asymptotic compactness allows us take some more general assumptions
than \eqref{1.11}-\eqref{1.12}.

For convenience, hereafter let $|\cdot|_p$ be the norm of
$L^p(\Omega)~(1\leqslant p < \infty)$, and $C$ a general positive
constant, which may be different in different estimates.

\medskip

This paper is organized as follows. We present some background
materials in \S 2, then prove a criterion on existence of pullback
attractors in \S 3, and a technical method  for verifying pullback
asymptotic compactness   is presented in \S 4. Finally, in \S 5,
these abstract results are applied to a non-autonomous wave
equation with nonlinear damping and critical nonlinearity, to
obtain the existence of pullback attractors. We conclude the paper
with some remarks in \S 6.

\section{Preliminaries}

\subsection{Kuratowski measure of non-compactness}

We briefly review the basic concept about the Kuratowski measure
of non-compactness and recall its basic properties, which will be
used to establish a criterion for the existence of pullback
attractors.
\begin{definition}(\cite{Ha,SY}) \label{d2.1}
Let $X$ be a complete metric space and $A$ be a bounded subset of
$X$. The Kuratowski measure of non-compactness $\kappa(A)$ of $A$
is defined as
\[
\kappa(A)=\inf\{\delta>0|\text{ A has a finite open cover of sets of
diameter} < \delta\}.
\]
\end{definition}

If $A$ is a nonempty, unbounded set in $X$, then we define
$\kappa(A)=\infty$.

The properties of $\kappa(A)$, which we will use in this paper, are
given in the following lemmas:
\begin{lemma}(\cite{Ha,SY}) \label{l2.2}
The Kuratowski measure of non-compactness $\kappa(A)$ on a
complete metric space $X$ satisfies the following properties:
\begin{enumerate}
\item[(1)]~$\kappa(A)=0$ if and only if $\bar{A}$ is compact,
where $\bar{A}$ is the closure of $A$;

\item[(2)]~$\kappa(\bar{A})=\kappa(A)$, $\kappa(A\cup B)=\max\{\kappa(A),\,\kappa(B)\}$;

\item[(3)]~If $A\subset B$, then $\kappa(A)\leqslant \kappa(B)$;

\item[(4)]~If $A_t$ is a family of nonempty, closed, bounded sets
defined for $t>r$ that satisfy $A_t \subset A_s$, whenever $s \leq
t$, and $\kappa(A_t)\to 0$, as $t \to \infty$, then
$\underset{t>r}{\cap}A_t$ is a nonempty, compact set in $X$.
\end{enumerate}
If in addition, $X$ is a Banach space, then the following estimate
is valid:
\begin{enumerate}
\item[(5)]~$\kappa(A+B)\leqslant \kappa(A)+\kappa(B)$\quad for any
bounded sets $A,B$ in $ X$.
\end{enumerate}
\end{lemma}

\subsection{Some useful properties for nonlinear damping function}

In the following, we will recall some simple properties of the
nonlinear damping function $h$, which will be used later.

\begin{lemma}(\cite{Fe,Kh})\label{l2.3}
Let $h$ satisfy \eqref{1.4} and \eqref{1.5}. Then for any $\delta
>0$, there exists a constant $C_{\delta}$   depending on $\delta$
such that
\[
|u-v|^2\leqslant \delta + C_{\delta}(h(u)-h(v))(u-v)\quad \text{for
all}~u,v\in \mathbb{R}.
\]
\end{lemma}

Moreover, condition \eqref{1.6} implies that
\[
|h(s)|^{\frac{1}{p}}\leqslant C(1+|s|).
\]
Therefore, we have
\[
|h(s)|^{\frac{p+1}{p}}=|h(s)|^{\frac{1}{p}}\cdot |h(s)|\leqslant
C(1+|s|)|h(s)|\leqslant C|h(s)|+Ch(s)\cdot s.
\]
Combining this estimate with the Young's inequality and
\eqref{1.4}, we further obtain that
\begin{equation}\label{2.1}
|h(s)|^{\frac{p+1}{p}}\leqslant C(1+h(s)\cdot s)\quad \text{for
all}~s\in \mathbb{R},
\end{equation}
where the constant $C$ is independent of $s$.

\section{Criterion for the existence of pullback attractors}

In this section, we first recall a few basic concepts  for
non-autonomous dynamical systems, including  pullback
$\kappa-$contraction, pullback asymptotic compactness and pullback
attractor. Then we present criteria for existence of pullback
attractors, in terms of $\kappa-$contraction or pullback
asymptotic compactness.

 Let $X$ be a complete metric space, which is the state space for
 a non-autonomous dynamical system (NDS).
 As in \cite{BCD,Ch,CKS}, we define a
non-autonomous dynamical system  in terms of a cocycle mapping
$\phi$: $\mathbb{R}^+\times \Sigma\times X \to X$ which is driven
by an \emph{autonomous} dynamical system $\theta$ acting on a
parameter space $\Sigma$. In details, $\theta=\{\theta_t\}_{t\in
\mathbb{R}}$ is a autonomous dynamical system on $\Sigma$, i.e., a
group of homeomorphisms under composition on $\Sigma$ with the
properties
\begin{enumerate}
\item[(i)] $\theta_0(\sigma) = \sigma$ for all $\sigma \in \Sigma$;

\item[(ii)] $\theta_{t+\tau}(\sigma) = \theta_t(\theta_{\tau}(\sigma))$ for
all $t,\tau\in \mathbb{R}$.
\end{enumerate}
The cocycle mapping $\phi$ satisfies
\begin{enumerate}
\item[(i)] $\phi(0, \sigma; x) = x$ for all $(\sigma,\, x)\in \Sigma \times X$;

\item[(ii)] $\phi(s + t,\, \sigma; x) = \phi(s, \theta_t(\sigma); \phi(t, \sigma;
x))$ for all $s, t\in \mathbb{R}^+$ and all $(\sigma, x)\in \Sigma
\times X$.
\end{enumerate}

Sometimes we say $\phi$ is a cocycle with respect to (w.r.t.)
$\theta$ and denote this by $(\phi, \theta)$.
 If,
in addition, the mapping $\phi(t, \sigma; \cdot):~X\to X$ is
continuous for each $\sigma \in \Sigma$ and $t\geqslant 0$, then
we call $\phi$ is a {\it continuous cocycle}. If the mapping
$\phi(t, \sigma; \cdot):~X\to X$ is \emph{norm-to-weak continuous}
for each $\sigma \in \Sigma$ and $t\geqslant 0$, that is, for each
$\sigma \in \Sigma$ and $t\geqslant 0$, norm convergence $x_n \to
x$ in $X$ implies weak convergence $\phi(t, \sigma;
x_n)\rightharpoonup \phi(t, \sigma; x)$, then we call $\phi$ is a
{\it norm-to-weak continuous cocycle}. A continuous cocycle is
obviously also a norm-to-weak continuous cocycle.

For convenience, hereafter, we will use the following notations:
\[
\mathcal{B}\overset{\vartriangle}{=}\{B~|~B~\text{is bounded
in}~X\}; \quad
\phi(t,\,\sigma;\,B)\overset{\vartriangle}{=}\{\phi(t,\,\sigma;\,x_0)~
|~x_0 \in B\}.
\]

\begin{definition}\label{d3.1} (\cite{Ch})
A family of bounded sets $\mathscr{B} = \{B_\sigma\}_{\sigma\in
\Sigma}$ of $X$ is called a bounded pullback absorbing set for the
cocycle $\phi$ with respect to (w.r.t.) $\theta$, if for any
$\sigma \in \Sigma$ and any $B\in \mathcal{B}$, there exists
$T=T(\sigma, B)\geqslant 0$ such that
\[
\phi(t,\, \theta_{-t}(\sigma); B)\subset B_{\sigma}\quad \text{for
all}~t\geqslant T.
\]
\end{definition}

\begin{definition}\label{d3.4}(\cite{Ch})
(\textbf{Pullback  attractor})\\
A family of nonempty compact sets $\mathscr{A} =
\{\mathcal{A}_\sigma\}_{\sigma\in \Sigma}$ of $X$ is called a
pullback  attractor for the cocycle $\phi$ w.r.t. $\theta$, if for
all $\sigma \in \Sigma$, it satisfies
\begin{enumerate}
\item[(i)] $\phi(t, \sigma; \mathcal{A}_{\sigma}) = \mathcal{A}_{\theta_t(\sigma)}$ for all $t \in
\mathbb{R}^+$ \quad ($\phi-$invariance);

\item[(ii)] $\lim\limits_{t\to +\infty} dist_X(\phi(t;
\theta_{-t}(\sigma); B),\,\mathcal{A}_{\sigma}) = 0$ for all bounded
set $B\subset X$.
\end{enumerate}
Often, $\mathcal{A}_\sigma$ is called a fiber at parameter $\sigma
\in \Sigma$.
\end{definition}

\begin{definition}\label{d3.5}(\cite{Ch})
Let $\phi$ be a cocycle w.r.t. $\theta$ on $\mathbb{R}^+\times
\Sigma \times X$, and let $B\in \mathcal{B}$. We define the
pullback $\omega$-limit set $\omega_{\sigma}(B)$ as   follows
\[
\omega_{\sigma}(B)=\bigcap_{s\geqslant
0}\overline{\bigcup_{t\geqslant s}\phi(t, \theta_{-t}(\sigma);
B)},  \;\;\;\;   \sigma \in \Sigma ,
\]
where $\overline{A}$ means the closure of $A$ in $X$.
\end{definition}

If the parameter space $\Sigma$  contains only one element
$\sigma_0$ and $\theta_t(\sigma_0)\equiv\sigma_0$ for all $t\in
\mathbb{R}$, then $\varphi$ reduces to a semigroup and all the
concepts in $Definitions$ \ref{d3.1}-\ref{d3.5} coincide with the
corresponding concepts   in autonomous systems. Especially, in the
autonomous case, the pullback attractor coincides with the global
attractor; see \cite{BV,Robinson,SY,Te, Ilyin}. Moreover,
Chepyzhov \& Vishik \cite{CV} define the concept of kernel
sections for non-autonomous dynamical systems, which correspond to
the fibers $\mathcal{A}_{\sigma}$ in the above Definition
\ref{d3.4} of a pullback attractor. Furthermore, similar to the
autonomous cases, we have also the following equivalent
characterization about the pullback $\omega$-limit set.

\begin{lemma}\label{l3.6} (\cite{Ch})
For any $B\subset \mathcal{B}$ and any $\sigma\in \Sigma$, $x_0\in
\omega_{\sigma}(B)$ if and only if there exist $\{x_n\}\subset B$
and $\{t_n\}\subset \mathbb{R}^+$ with $t_n \to +\infty$ as $n\to
\infty$, such that
\[
\phi(t_n, \theta_{-t_n}(\sigma); x_n) \to x_0\quad \text{as}~n\to
\infty.
\]
\end{lemma}

Now  we define the pullback $\kappa$-contracting cocycle in terms
 of the Kuratowski non-compactness measure:
\begin{definition}\label{d3.7} (\textbf{$\kappa$-contracting cocycle})\\
Let $\phi$ be a cocycle w.r.t. $\theta$ on $\mathbb{R}^+\times
\Sigma \times X$. Then $\phi$ is called pullback
$\kappa$-contracting if for any $\varepsilon>0$, $\sigma\in
\Sigma$ and any $B\in \mathcal{B}$, there is a $T=T(\varepsilon,
\sigma, B)\geqslant 0$ such that
\[
\kappa_{_X}\left(\phi(t, \theta_{-t}(\sigma); B)\right) \leqslant
\varepsilon \quad \text{for all}~t\geqslant T.
\]
\end{definition}

From the definitions above, we have the following basic fact.
\begin{lemma}\label{c3.8}
Let $\phi$ be a cocycle w.r.t. $\theta$ on $\mathbb{R}^+\times
\Sigma \times X$. If $\phi$ has a pullback attractor, then $\phi$
has a bounded pullback absorbing set and $\phi$ is pullback
$\kappa$-contracting.
\end{lemma}

We introduce another definition, needed for characterizations of
existence of pullback attractors later.
\begin{definition}\label{d3.2}
(\textbf{Nested   pullback absorbing set})\\
A family of bounded sets $\mathscr{B} = \{B_\sigma\}_{\sigma\in
\Sigma}$ of $X$ is called a nested bounded pullback absorbing set
for $\phi$ w.r.t. $\theta$ if $\mathscr{B}$ is a bounded pullback
absorbing set, and, moreover, $B_{\sigma}$ satisfy the nested
relation: $B_{\theta_{-t}(\sigma)}\subset B_{\sigma}$ for any
$t\geqslant 0$ and any $\sigma \in \Sigma$.
\end{definition}

\begin{remark}
This nested relation appears in some systems arising in physical
applications.
 For example, the non-autonomous systems considered in \cite{CKR,Ch,
CV}  have   nested bounded pullback absorbing sets.

\end{remark}

In the following, we will present some characterizations for the
pullback $\kappa$-contracting cocycles.
\begin{lemma}\label{l3.9}
Let $\phi$ be a $\kappa$-contracting cocycle w.r.t. $\theta$ on
$\mathbb{R}^+\times \Sigma \times X$ and have a nested bounded
pullback absorbing set $\mathscr{B}=\{B_{\sigma}\}_{\sigma \in
\Sigma}$. Then for every $\sigma\in \Sigma$,   every bounded
sequence $\{x_n\}_{n=1}^{\infty}\subset X$ and every time sequence
$\{t_n\}\subset \mathbb{R}^+$ with $t_n\to +\infty$ as $n \to
\infty$, we have
\begin{enumerate}
\item[(i)] $\{\phi(t_n,\,\theta_{-t_n}(\sigma);\,
x_n)\}_{n=1}^{\infty}$ is pre-compact in $X$;

\item[(ii)] all clusters of $\{\phi(t_n,\,\theta_{-t_n}(\sigma);\,
x_n)\}_{n=1}^{\infty}$ are contained in
$\omega_{\sigma}(B_{\sigma})$, that is, if
\begin{equation*}
\phi(t_{n_j},\,\theta_{-t_{n_j}}(\sigma);\, x_{n_j}) \to x_0 \quad
\text{as}~j \to \infty,
\end{equation*}
then $x_0\in \omega_{\sigma}(B_{\sigma})$;

\item[(iii)] $\omega_{\sigma}(B_{\sigma})$ is nonempty and compact in $X$.
\end{enumerate}
\end{lemma}

{\bf Proof.} $(i).$ Denote $\{x_n\}_{n=1}^{\infty}$ by $B$. For
any $\varepsilon>0$ and for each $\sigma \in \Sigma$, by the
definition of pullback $\kappa$-contracting cocycle, we know that
there exists a $T_0=T_0(\varepsilon,\,\sigma,B_{\sigma}) >0$ such
that
\begin{equation}\label{3.1}
\kappa_{_X}\left(\phi(t, \theta_{-t}(\sigma); B_{\sigma})\right)
\leqslant \varepsilon \quad \text{for all}~t\geqslant T_0
\end{equation}
and there exists also a $T_1=T_1(\varepsilon, \, \sigma,B)$ such
that
\begin{equation}\label{3.2}
\phi(t+T_1,\,\theta_{-(t+T_1)}(\theta_{-T_0}(\sigma));\, B) \subset
B_{\theta_{-T_0}(\sigma)} \subset B_{\sigma}\quad \text{for all}~
t\geqslant 0.
\end{equation}
Hence, for any $t\geqslant 0$, we have
\begin{align}\label{3.3}
\phi(t+T_1+T_0,\,& \theta_{-(t+T_1+T_0)}(\sigma);\, B)
\nonumber \\
& = \phi(T_0,\,\theta_{-T_0}(\sigma);\,
\phi(t+T_1,\,\theta_{-(t+T_1)}(\theta_{-T_0}(\sigma)); B))
\nonumber \\
& \subset \phi(T_0,\,\theta_{-T_0}(\sigma);\,
B_{\theta_{-T_0}(\sigma)}) \nonumber \\
& \subset \phi(T_0,\,\theta_{-T_0}(\sigma);\, B_{\sigma}),
\end{align}
and then
\begin{equation}\label{3.4}
\bigcup_{t\geqslant T_0+T_1}\phi(t, \theta_{-t}(\sigma);\, B)
\subset \phi(T_0,\,\theta_{-T_0}(\sigma);\, B_{\sigma}).
\end{equation}

Therefore, combining \eqref{3.1} and \eqref{3.4}, we have
\begin{equation}\label{3.5}
\kappa_{_X}\left(\bigcup_{t\geqslant T_0+T_1}\phi(t,
\theta_{-t}(\sigma);\, B)\right) \leqslant \varepsilon.
\end{equation}

Then by the properties $(1),(2)$ of $Lemma$ \ref{l2.2} and
$\{\phi(t_n,\,\theta_{-t_n}(\sigma);\, x_n)\}_{n=n_0}^{\infty}
\subset \bigcup_{t\geqslant T_0+T_1}\phi(t, \theta_{-t}(\sigma);\,
B)$ for some $n_0$, we know that
$\kappa_X(\{\phi(t_n,\,\theta_{-t_n}(\sigma);\,
x_n)\}_{n=1}^{\infty}) \leqslant \varepsilon$. Hence by the
arbitrariness of $\varepsilon$ and property $(1)$ of $Lemma$
\ref{l2.2}, we conclude that
$\{\phi(t_n,\,\theta_{-t_n}(\sigma);\, x_n)\}_{n=1}^{\infty}$ is
pre-compact in $X$.

$(ii).$ Let $x_0$ be a cluster of
$\{\phi(t_n,\,\theta_{-t_n}(\sigma);\, x_n)\}_{n=1}^{\infty}$, we
need to show that $x_0\in \omega_{\sigma}(B_{\sigma})$. Without
loss of generality, we assume that
$\phi(t_n,\,\theta_{-t_n}(\sigma);\, x_n) \to x_0$ as $n\to
\infty$.

We claim first that for each sequence $\{s_m\}_{m=1}^{\infty}
\subset \mathbb{R}^+$ satisfying $s_m \to \infty$ as $m\to \infty$,
we can find two sequences $\{t_{n_m}\}_{m=1}^{\infty}\subset
\{t_{n}\}_{n=1}^{\infty}$ and $\{y_m\}_{m=1}^{\infty}\subset
B_{\sigma}$ satisfying $t_{n_m} \to \infty$ as $m\to \infty$, such
that
\begin{equation}\label{3.6}
\phi(s_m,\,\theta_{-s_m}(\sigma);\,y_m)=\phi(t_{n_m},\,\theta_{-t_{n_m}}(\sigma);\,x_{n_m}).
\end{equation}
Indeed, for each $m\in \mathbb{N}$, we can take $n_m$ so large that
$t_{n_m} \geqslant s_m$ and
\[
y_m\overset{\vartriangle}{=} \phi(t_{n_m}-s_m,\,
\theta_{-(t_{n_m}-s_m)} (\theta_{-s_m}(\sigma)); \,x_{n_m}) \in
B_{\theta_{-s_m}(\sigma)} \subset B_{\sigma}.
\]
Therefore,
\begin{align}\label{3.7}
\phi(t_{n_m},\,& \theta_{-t_{n_m}}(\sigma);\, x_{n_m}) \nonumber \\
& = \phi(s_m+(t_{n_m}-s_m),\,
\theta_{-(s_m+(t_{n_m}-s_m))}(\sigma);\, x_{n_m})
\nonumber \\
& = \phi(s_m,\,\theta_{-s_m}(\sigma);\,\phi(t_{n_m}-s_m,\,
\theta_{-(t_{n_m}-s_m)}(\theta_{-s_m}(\sigma));\,x_{n_m})) \nonumber \\
& = \phi(s_m,\,\theta_{-s_m}(\sigma);\,y_m).
\end{align}

Hence,
\[
\lim\limits_{m\to \infty} \phi(s_m,\,\theta_{-s_m}(\sigma);\,y_m) =
\lim\limits_{m\to \infty}\phi(t_{n_m},\,\theta_{-t_{n_m}}
(\sigma);\,x_{n_m}) =x_0,
\]
and $y_m\in B_{\sigma}$ for each $m\in \mathbb{N}$, which implies,
by the definition of $\omega_{\sigma}(B_{\sigma})$, that $x_0 \in
\omega_{\sigma}(B_{\sigma})$.

$(iii).$ The fact that    $\omega_{\sigma}(B_{\sigma})$ is
nonempty is obvious. Substitute $B$ by $B_{\sigma}$ in
\eqref{3.2}-\eqref{3.5}, we   obtain that there exists a
$T_2=T_2(\varepsilon,B_{\sigma},\sigma)$ such that
\begin{equation*}
\kappa_{_X}\left(\overline{\bigcup_{t\geqslant T_0+T_2}\phi(t,
\theta_{-t}(\sigma);\, B_{\sigma})}\right) =
\kappa_{_X}\left(\bigcup_{t\geqslant T_0+T_2}\phi(t,
\theta_{-t}(\sigma);\, B_{\sigma})\right) \leqslant \varepsilon.
\end{equation*}
Then by the definition of pullback $\omega$-limit set and property
$(4)$ of $Lemma$ \ref{l2.2}, we know that
$\omega_{\sigma}(B_{\sigma})$ is compact in $X$. $\hfill
\blacksquare$

A criterion for the existence of pullback attractors is then
obtained by means of $\kappa$-contraction.

\begin{theorem}\label{t3.10}
(\textbf{Sufficient condition for existence of pullback attractors})\\
 Let
$\phi$ be a continuous cocycle w.r.t. $\theta$ on
$\mathbb{R}^+\times \Sigma \times X$. Then $(\phi, \theta)$ has a
pullback attractor provided that
\begin{enumerate}
\item[(i)]\ $(\phi, \theta)$ has a nested bounded pullback absorbing set $\mathscr{B}=\{B_{\sigma}\}_{\sigma \in \Sigma}$;

\item[(ii)]\ $(\phi, \theta)$ is pullback $\kappa$-contracting.
\end{enumerate}
\end{theorem}
{\bf Proof.} For any $\sigma \in \Sigma$, we consider a family of
 $\omega$-limit sets $\mathscr{B}=\{B_{\sigma}\}_{\sigma\in
\Sigma}$:
\begin{equation*}
\omega_{\sigma}(B_{\sigma})=\bigcap_{s\geq
0}\overline{\bigcup_{t\geq  s}\phi(t, \theta_{-t}(\sigma);
B_{\sigma})}, \quad \sigma \in \Sigma.
\end{equation*}

By $Lemma$ \ref{l3.9} we know that $\omega_{\sigma}(B_{\sigma})$ is
nonempty and compact in $X$ for each $\sigma \in \Sigma$.

In the following, we will prove that
$\mathscr{A}=\{\omega_{\sigma}(B_{\sigma})\}_{\sigma \in \Sigma}$
is a pullback attractor of $(\phi,\,\theta)$, which will be
accomplished in two steps.

{\it Claim 1. For each $\sigma\in \Sigma$ and any $B\in
\mathcal{B}$, we have
\begin{equation*}
\lim\limits_{t\to +\infty} dist_X(\phi(t,\, \theta_{-t}(\sigma);
B),\,\omega_{\sigma}(B_{\sigma})) = 0.
\end{equation*}
}

In fact, if {\it Claim 1} is not true, then there exist
$\varepsilon_0>0$, $\{x_n\}_{n=1}^{\infty}\subset B$ and $\{t_n\}$
with $t_n\to +\infty$ as $n\to \infty$, such that
\begin{equation}\label{3.8}
dist_X(\phi(t_n,\, \theta_{-t_n}(\sigma);\,
x_n),\,\omega_{\sigma}(B_{\sigma})) \geqslant \varepsilon_0 \quad
\text{for}~n=1,2,\cdots.
\end{equation}

However, thanks to $Lemma$ \ref{l3.9}, we know that
$\{\phi(t_n,\,\theta_{-t_n}(\sigma);\,x_n)\}_{n=1}^{\infty}$ is
pre-compact in $X$. Without loss of generality, we assume that
\begin{equation}\label{3.9}
\phi(t_n,\,\theta_{-t_n}(\sigma);\,x_n) \to x_0  \quad \text{as}~n
\to \infty.
\end{equation}
Then $x_0 \in \omega_{\sigma}(B_{\sigma})$, which is a contraction
with \eqref{3.8}. This complete the proof of {\it Claim 1}.

{\it Claim 2. $\mathscr{A}=\{\omega_{\sigma}(B_{\sigma})\}_{\sigma
\in \Sigma}$ is $\phi$ invariant, that is,
\begin{equation*}
\phi(t,\,\sigma;\, \omega_{\sigma}(B_{\sigma}))=
\omega_{\theta_t(\sigma)}(B_{\theta_t(\sigma)})\quad \text{for
all}~t\geqslant 0,~\sigma\in \Sigma.
\end{equation*}
}

We first take $x \in \phi(t,\,\sigma;\,
\omega_{\sigma}(B_{\sigma}))$.

Then there is a $y\in \omega_{\sigma}(B_{\sigma})$ such that
$x=\phi(t,\,\sigma;\,y)$, and by the definition of $y$, there exist
$\{y_n\}\subset B_{\sigma} \subset B_{\theta_t(\sigma)}$ and $t_n$
with $t_n \to \infty$ as $n \to \infty$ such that
$y=\lim\limits_{n\to
\infty}\phi(t_n,\,\theta_{-t_n}(\sigma);\,y_n)$.

Therefore, by the continuity of $\phi$, as $n \to \infty$,
\begin{equation}\label{3.10}
\phi(t_n+t,\,\theta_{-(t_n+t)}(\theta_t(\sigma));\,y_n) =
\phi(t,\,\sigma;\,\phi(t_n,\,\theta_{-t_n}(\sigma);\,y_n)) \to
\phi(t,\,\sigma;\,y)=x.
\end{equation}

On the other hand, from $Lemma$ \ref{l3.9}, we know that
$\{\phi(t_n+t,\,\theta_{-(t_n+t)}(\theta_t(\sigma));\,y_n)\}_{n=1}^{\infty}$
is pre-compact in $X$. Without loss of generality, we assume that
\begin{equation*}
\phi(t_n+t,\,\theta_{-(t_n+t)}(\theta_t(\sigma));\,y_n) \to x_0
\in \omega_{\theta_t(\sigma)}(B_{\theta_t(\sigma)}) \quad
\text{as} ~n\to \infty.
\end{equation*}
Then by the uniqueness of limitation, we have $x=x_0$, which
implies that $x\in
\omega_{\theta_t(\sigma)}(B_{\theta_t(\sigma)})$, and thus
\begin{equation}\label{3.11}
\phi(t,\,\sigma;\,\omega_{\sigma}(B_{\sigma})) \subset
\omega_{\theta_t(\sigma)}(B_{\theta_t(\sigma)}).
\end{equation}
Now, we only need to prove the converse inclusion relation.

Let $z\in \omega_{\theta_t(\sigma)}(B_{\theta_t(\sigma)})$. Then
there exist $\{z_n\}\subset B_{\theta_t(\sigma)}$ and $t_n$ with
$t_n \to \infty$ as $n \to \infty$ such that $z=\lim\limits_{n\to
\infty}\phi(t_n,\,\theta_{-t_n}(\theta_t(\sigma));\,z_n)$.

Since $\{z_n\}\subset B_{\theta_t(\sigma)}$ is bounded, from
$Lemma$ \ref{l3.9}, we know that
$\{\phi(t_n-t,\,\theta_{-(t_n-t)}(\sigma);\,z_n)\}_{n=1}^{\infty}$
is pre-compact in $X$. Without loss of generality, we assume that
$\phi(t_n-t,\,\theta_{-(t_n-t)}(\sigma);\,z_n) \to x_0 \in
\omega_{\sigma}(B_{\sigma})$ as $n \to \infty$. Then by the
continuity of $\phi$, we have
\begin{align}\label{3.12}
\phi(t,\,\sigma;\,x_0)\leftarrow
\phi(t,&\,\sigma;\,\phi(t_n-t,\,\theta_{-(t_n-t)}(\sigma);\,z_n))
\nonumber \\
& =\phi(t_n, \, \theta_{-(t_n-t)}(\sigma);\,z_n) \nonumber \\
& = \phi(t_n, \, \theta_{-(t_n)}(\theta_t(\sigma));\,z_n) \to z.
\end{align}
Hence, $z= \phi(t,\,\sigma;\,x_0)$ with $x_0\in
\omega_{\sigma}(B_{\sigma})$, which implies
\begin{equation}\label{3.13}
\omega_{\theta_t(\sigma)}(B_{\theta_t(\sigma)}) \subset
\phi(t,\,\sigma;\,\omega_{\sigma}(B_{\sigma})).
\end{equation}

Combining \eqref{3.11} and \eqref{3.13} we know that {\it Claim 2}
is true.

From {\it Claim 1} and {\it Claim 2}, we complete the proof of
$Theorem$ \ref{t3.10}. $\hfill \blacksquare$

\begin{remark}\label{r3.11}
In the proof of $Theorem$ \ref{t3.10}, the continuity of the cocycle
$\phi(t,\,\sigma;\,\cdot):\, X\to X$ can be replaced by the
``weaker" continuity; see \eqref{3.10} and \eqref{3.12}. That is,
the above proof holds for norm-to-weak continuous cocycles; see
\cite{ZYS} for autonomous cases.
\end{remark}

Similar to the definition in Caraballo et al \cite{CLR}, we define
the following {\it pullback asymptotic compactness} for NDS.
\begin{definition}\label{d3.12}
Let $\phi$ be a cocycle w.r.t. $\theta$ on $\mathbb{R}^+\times
\Sigma \times X$. Then $\phi$ is called pullback asymptotically
compact, if for each $\sigma\in \Sigma$,     every bounded
sequence $\{x_n\}_{n=1}^{\infty}$, and every time sequence
$\{t_n\}\subset \mathbb{R}^+$ with $t_n\to +\infty$ as $n \to
\infty$, $\{\phi(t_n,\, \theta_{-t_n}(\sigma); \,
x_n)\}_{n=1}^{\infty}$ is pre-compact in $X$.
\end{definition}

In the framework of pullback attractors, the pullback asymptotic
compactness may not be equivalent to the $\kappa$-contraction if
the cocycle only has a general bounded pullback absorbing set. In
fact, we need this bounded pullback absorbing set to satisfy an
additional nesting condition; see the next theorem.

From $Lemma$ \ref{l3.9} we know that if $\phi$ has a nested
bounded pullback absorbing set, then $\phi$ being pullback
$\kappa$-contracting implies that $\phi$ being pullback
asymptotically compact; furthermore, in the proof of $Theorem$
\ref{t3.10}, we note that we indeed only used the pullback
asymptotic compactness. This, combining with $Lemma$ \ref{c3.8},
implies the following criterion.

\begin{theorem} \label{t3.13}(\textbf{Criterion for existence of pullback
attractor})\\ Let $\phi$ be a norm-to-weak continuous cocycle
w.r.t. $\theta$ on $\mathbb{R}^+\times \Sigma \times X$ such that
$(\phi, \theta)$ has a nested bounded pullback absorbing set. Then
$(\phi, \theta)$ has a pullback attractor if and only if $(\phi,
\theta)$ is pullback $\kappa$-contracting, or equivalently,
$(\phi, \theta)$ is pullback asymptotically compact.
\end{theorem}
That is, under the assumption that $(\phi, \theta)$ has a nested
bounded pullback absorbing set, pullback $\kappa$-contraction is
equivalent to pullback asymptotic compactness.

On the other hand, the authors in \cite{CLR} have proven that
$(\phi, \theta)$ has a pullback attractor provided that $(\phi,
\theta)$ is pullback asymptotically compact and has a bounded
pullback absorbing set (see $Theorem$ 7 of \cite{CLR}). In fact,
from the definition of pullback attractor, $Lemma$ \ref{c3.8},
Theorems \ref{t3.10} and \ref{t3.13},  we   observe that these
conditions are also necessary. We summarize this result in the
following theorem.

\begin{theorem}\label{t3.14}
(\textbf{Another criterion for existence of pullback attractor})\\
Let $\phi$ be a norm-to-weak continuous cocycle w.r.t. $\theta$ on
$\mathbb{R}^+\times \Sigma \times X$. Then $(\phi, \theta)$ has a
pullback attractor if and only if $(\phi, \theta)$ is pullback
asymptotically compact and has a bounded pullback absorbing set.
\end{theorem}

$Theorem$ \ref{t3.13} and $Theorem$ \ref{t3.14} show that pullback
asymptotically compact is stronger than pullback
$\kappa-$contracting to some extent, which is different from the
autonomous cases. Note also that $Theorem$ 3.14 is a slight
improvement of $Theorem$ 7 in \cite{CLR}, from continuous cocycles
to   ``weakly" continuous cocycles (i.e., norm-to-weak continuous
cocycles).

Although we only use the pullback asymptotic  compactness in our
later applications in \S 5, we think that the pullback
$\kappa-$contraction criterion for existence of pullback
attractors for ``weakly" continuous cocycles (i.e.,
 norm-to-weak continuous cocycles) is of independent interest and will be
useful for other non-autonomous dynamical systems. Another reason
to present the $\kappa-$contraction criterion here is that we like
to highlight a difference with the \emph{autonomous} systems: In
non-autonomous systems, the pullback asymptotic  compactness
criterion  and the pullback $\kappa-$contraction criterion, for
existence of pullback attractors, are not equivalent unless when
there exists a \emph{nested} bounded absorbing set ($Theorem$
\ref{t3.13}).


\medskip

We also remark that     the definitions and results in this
section can be expressed in the framework of \emph{processes},
instead of cocycles, as   in   \cite{CV}.


\section{A technical method for verifying pullback asymptotic compactness}

We now present a convenient method for verifying the pullback
asymptotic compactness for the cocycle generated by
\emph{non-autonomous} hyperbolic type of equations, in order to
apply $Theorem$ \ref{t3.13} to obtain existence of pullback
attractors in the next section. This method is partially motivated
by the methods in \cite{CL2,CL3,Kh} in some sense; see also in
\cite{SYZ2}. In \cite{CL3}, the authors present a general abstract
framework for asymptotic dynamics of \emph{autonomous} wave
equations.

\begin{definition}(\cite{SYZ2})
 Let $X$ be a Banach space  and $B$ be a bounded subset of $X$.
We call a function $\psi(\cdot,\cdot)$, defined on $X\times X$, a
contractive function on $B\times B$ if for any sequence
$\{x_n\}_{n=1}^{\infty}\subset B$, there is a subsequence
$\{x_{n_k}\}_{k=1}^{\infty} \subset \{x_n\}_{n=1}^{\infty}$ such
that
\begin{equation*}
\lim_{k\to \infty}\lim_{l\to \infty}\psi(x_{n_k},\,x_{n_l})=0 .
\end{equation*}
We denote the set of all   contractive functions on $B\times B$ by
$Contr(B)$.
\end{definition}

\begin{theorem}\label{t4.2}
\textbf{(Technique for verifying pullback asymptotic compactness)}
\\
Let $\phi$ be a cocycle w.r.t. $\theta$ on $\mathbb{R}^+\times
\Sigma \times X$ and have a nested bounded pullback absorbing set
$\mathscr{B}=\{B_{\sigma}\}_{\sigma \in \Sigma}$. Moreover, assume
that for any $\varepsilon >0$ and each $\sigma\in \Sigma$, there
exist $T=T(B_{\sigma}, \varepsilon)$ and
$\psi_{T,\,\sigma}(\cdot,\cdot) \in Contr(B_{\sigma})$ such that
\begin{equation*}
\|\phi(T,\,\theta_{-T}(\sigma);x)-\phi(T,\,\theta_{-T}(\sigma);y)\|\leqslant
\varepsilon + \psi_{T,\,\sigma}(x,\,y) \quad \text{for all}~x,y\in
B_{\sigma},
\end{equation*}
where  $\psi_{T,\,\sigma}$ depends on $T$ and $\sigma$. Then $\phi$
is pullback asymptotically compact in $X$.
\end{theorem}

{\bf Proof.} Let $\{y_n\}_{n=1}^{\infty}$ be a bounded sequence of
$X$ and $\{t_n\}\subset \mathbb{R}^+$ with $t_n \to \infty$ as
$n\to \infty$}. We need to show that
\begin{align}\label{4.1}
\text{$\{\phi(t_n,\,\theta_{-t_n}(\sigma);\,y_n)\}_{n=1}^{\infty}$
is precompact in $X$ for each $\sigma\in \Sigma$}.
\end{align}

In the following, we will prove that
$\{\phi(t_n,\,\theta_{-t_n}(\sigma);\,y_n)\}_{n=1}^{\infty}$ has a
convergent subsequence via diagonal methods (e.g., see \cite{Kh}).

Taking $\varepsilon_m>0$ with $\varepsilon_m \to 0$ as $m\to
\infty$.

At first, for $\varepsilon_1$, by the assumptions, there exist
$T_1=T_1(\varepsilon_1)$ and $\psi_{1}(\cdot,\cdot) \in
Contr(B_{\sigma})$ such that
\begin{equation}\label{4.2}
\|\phi(T_1,\,\theta_{-T_1}(\sigma);\,x)-\phi(T_1,\,\theta_{-T_1}(\sigma);\,y)\|\leqslant
\varepsilon_1 + \psi_1(x,\,y) \quad \text{for all}~x,y\in
B_{\sigma},
\end{equation}
where $\psi_1$ depends on $T_1$ and $\sigma$.

Since $t_n\to \infty$, for such fixed $T_1$, without loss of
generality, we assume that $t_n$ is so large that
\begin{equation}\label{4.3}
\phi(t_n-T_1,\,\theta_{-(t_n-T_1)}(\theta_{-T_1}(\sigma));\,y_n) \in
B_{\theta_{-T_1}(\sigma)} \subset B_{\sigma}\quad \text{for
each}~n=1,2,\cdots.
\end{equation}
Set
$x_n=\phi(t_n-T_1,\,\theta_{-(t_n-T_1)}(\theta_{-T_1}(\sigma));\,y_n)$.
Then from \eqref{4.2} we have
\begin{align}\label{4.4}
\|\phi(t_n,\,\theta_{-t_n}&(\sigma);\,y_n)-\phi(t_m,\,\theta_{-t_m}(\sigma);\,y_m)\|
\nonumber \\
&= \|\phi(T_1,\theta_{-T_1}(\sigma); x_n)
-\phi(T_1,\theta_{-T_1}(\sigma);x_m)\| \leqslant \varepsilon_1 +
\psi_1(x_n,\,x_m).
\end{align}

Due to the definition of $Contr(B_{\sigma})$ and
$\psi_1(\cdot,\cdot) \in Contr(B_{\sigma})$, we know that
$\{x_n\}_{n=1}^{\infty}$ has a subsequence
$\{x^{(1)}_{n_k}\}_{k=1}^{\infty}$ such that
\begin{equation}\label{4.5}
\lim_{k\to \infty}\lim_{l\to
\infty}\psi_1(x^{(1)}_{n_k},\,x^{(1)}_{n_l}) \leqslant
\frac{\varepsilon_1}{2},
\end{equation}
and similar to  \cite{Kh}, we have
\begin{align*}
\lim_{k\to \infty}\sup_{p\in \mathbb{N}}&
\|\phi(t^{(1)}_{n_{k+p}},\,\theta_{-t^{(1)}_{n_{k+p}}}(\sigma);\,y^{(1)}_{n_{k+p}})
-\phi(t^{(1)}_{n_{k}},\,\theta_{-t^{(1)}_{n_{k}}}(\sigma);\,y^{(1)}_{n_{k}})\|
\nonumber \\
&\leqslant \lim_{k\to \infty}\sup_{p\in \mathbb{N}}\limsup_{l\to
\infty}
\|\phi(t^{(1)}_{n_{k+p}},\,\theta_{-t^{(1)}_{n_{k+p}}}(\sigma);\,y^{(1)}_{n_{k+p}})
-\phi(t^{(1)}_{n_{l}},\,\theta_{-t^{(1)}_{n_{l}}}(\sigma);\,y^{(1)}_{n_{l}})\|
\nonumber \\
& \qquad + \limsup_{k\to \infty}\limsup_{l\to \infty}
\|\phi(t^{(1)}_{n_{k}},\,\theta_{-t^{(1)}_{n_{k}}}(\sigma);\,y^{(1)}_{n_{k}})-\phi(t^{(1)}_{n_{l}},\,\theta_{-t^{(1)}_{n_{l}}}(\sigma);\,y^{(1)}_{n_{l}})\|
\nonumber \\
&\leqslant \varepsilon_1+\lim_{k\to \infty}\sup_{p\in
\mathbb{N}}\lim_{l\to
\infty}\psi_1(x^{(1)}_{n_{k+p}},\,x^{(1)}_{n_l}) +\varepsilon_1+
\lim_{k\to \infty}\lim_{l\to
\infty}\psi_1(x^{(1)}_{n_k},\,x^{(1)}_{n_l}),
\end{align*}
which, combining with \eqref{4.4} and \eqref{4.5}, implies that
\begin{align*}
\lim_{k\to \infty}\sup_{p\in \mathbb{N}}
\|\phi(t^{(1)}_{n_{k+p}},\,\theta_{-t^{(1)}_{n_{k+p}}}(\sigma);\,y^{(1)}_{n_{k+p}})
-\phi(t^{(1)}_{n_{k}},\,\theta_{-t^{(1)}_{n_{k}}}(\sigma);\,y^{(1)}_{n_{k}})\|
\leqslant 4\varepsilon_1.
\end{align*}

Therefore, there is a $K_1$ such that
\begin{align*}
\|\phi(t^{(1)}_{n_{k}},\,\theta_{-t^{(1)}_{n_{k}}}(\sigma);\,y^{(1)}_{n_{k}})
-\phi(t^{(1)}_{n_{l}},\,\theta_{-t^{(1)}_{n_{l}}}(\sigma);\,y^{(1)}_{n_{l}})\|
\leqslant 5\varepsilon_1 \quad \text{for all}~k,l\geqslant K_1.
\end{align*}

By induction, we obtain that, for each $m\geqslant 1$, there is a
subsequence \linebreak
$\{\phi(t^{(m+1)}_{n_{k}},\,\theta_{-t^{(m+1)}_{n_{k}}}(\sigma);\,y^{(m+1)}_{n_{k}})\}_{k=1}^{\infty}$
of
$\{\phi(t^{(m)}_{n_{k}},\,\theta_{-t^{(m)}_{n_{k}}}(\sigma);\,y^{(m)}_{n_{k}})\}_{k=1}^{\infty}$
and certain $K_{m+1}$ such that
\begin{equation*}
\|\phi(t^{(m+1)}_{n_{k}},\,\theta_{-t^{(m+1)}_{n_{k}}}(\sigma);\,y^{(m+1)}_{n_{k}})-\phi(t^{(m+1)}_{n_{l}},\,\theta_{-t^{(m+1)}_{n_{l}}}(\sigma);\,y^{(m+1)}_{n_{l}})\|
\leqslant 5\varepsilon_{m+1} \quad \text{for all}~k,l\geqslant
K_{m+1}.
\end{equation*}

Now, we consider the diagonal subsequence
$\{\phi(t^{(k)}_{n_{k}},\,\theta_{-t^{(k)}_{n_{k}}}(\sigma);\,y^{(k)}_{n_{k}})\}_{k=1}^{\infty}$.
Since for each $m\in \mathbb{N}$,
$\{\phi(t^{(k)}_{n_{k}},\,\theta_{-t^{(k)}_{n_{k}}}(\sigma);\,y^{(k)}_{n_{k}})\}_{k=m}^{\infty}$
is a subsequence of
$\{\phi(t^{(m)}_{n_{k}},\,\theta_{-t^{(m)}_{n_{k}}}(\sigma);\,y^{(m)}_{n_{k}})\}_{k=1}^{\infty}$,
then,
\begin{equation*}
\|\phi(t^{(k)}_{n_{k}},\,\theta_{-t^{(k)}_{n_{k}}}(\sigma);\,y^{(k)}_{n_{k}})-
\phi(t^{(l)}_{n_{l}},\,\theta_{-t^{(l)}_{n_{l}}}(\sigma);\,y^{(l)}_{n_{l}})\|
\leqslant 5\varepsilon_{m} \quad \text{for all}~k,l\geqslant
\max\{m,K_{m}\},
\end{equation*}
which, combining with $\varepsilon_m \to 0$ as $m\to \infty$,
implies that
$\{\phi(t^{(k)}_{n_{k}},\,\theta_{-t^{(k)}_{n_{k}}}(\sigma);\,y^{(k)}_{n_{k}})\}_{k=1}^{\infty}$
is a Cauchy sequence in $X$. This shows that
$\{\phi(t_n,\,\theta_{-t_n}(\sigma);\,y_n)\}_{n=1}^{\infty}$ is
precompact in $X$ for each $\sigma\in \Sigma$. We thus complete
the proof.
 $\hfill \blacksquare$ \\

Note that the nested properties and the contractive properties are
only used in \eqref{4.3} and \eqref{4.2} respectively. We have a
similar corollary for the cocycle without nested pullback absorbing
set, the proof is similar to that for $Theorem$ \ref{t4.2} above.
\begin{corollary}\label{c4.3}
Let $\phi$ be a cocycle w.r.t. $\theta$ on $\mathbb{R}^+\times
\Sigma \times X$ and have a bounded pullback absorbing set
$\mathscr{B}=\{B_{\sigma}\}_{\sigma \in \Sigma}$. Moreover, assume
that for any $\varepsilon >0$ and each $\sigma\in \Sigma$, there
exist $T=T(\sigma, \varepsilon)$ and $\psi_{T,\,\sigma}(\cdot,\cdot)
\in Contr(B_{\theta_{-T}(\sigma)})$ such that
\begin{equation*}
\|\phi(T,\,\theta_{-T}(\sigma);x)-\phi(T,\,\theta_{-T}(\sigma);y)\|\leqslant
\varepsilon + \psi_{T,\,\sigma}(x,\,y) \quad \text{for all}~x,y\in
B_{\theta_{-T}(\sigma)},
\end{equation*}
where $\psi_{T,\,\sigma}$ depends on $T$ and $\sigma$. Then $\phi$
is pullback asymptotically compact in $X$.
\end{corollary}

\section{Pullback attractors for a non-autonomous wave equation}

In this section, we prove the existence of the pullback attractor
for the non-autonomous wave system \eqref{1.1}-\eqref{1.3}, by
applying $Theorems$ \ref{t3.13} and \ref{t3.14}. We use the method
(via contractive functions) in \S 4 to verify the pullback
asymptotic compactness. This method appears to be very efficient for
non-autonomous wave or hyperbolic equations, while the approach in
\cite{CLR}, which is an energy method and is different from ours, is
very appropriate for some non-autonomous parabolic equations or wave
equations with \emph{linear} damping, e.g., see \cite{MRW}. In fact,
the approach in \cite{CLR} is an energy method and may be seen as a
non-autonomous generalization of  Ball's method \cite{Ba}.

\subsection{Mathematic setting}

We consider the non-autonomous wave system \eqref{1.1}-\eqref{1.3}
on the \emph{state space} $X=H_0^1(\Omega)\times L^2(\Omega)$.
 For each $g_0\in L_{loc}^2(\mathbb{R};\,L^2(\Omega))$,
we denote $\{g_0(s+t)|t\in \mathbb{R}\}$ by $\mathcal{H}_1(g_0)$.
For $f_0(v,s)$ satisfying \eqref{1.8}-\eqref{1.9}, we similarly
denote $\mathcal{H}_2(f_0)=\{f_0(\cdot,s+t)|t\in \mathbb{R}\}$.

Let $\Sigma=\mathcal{H}_2(f_0)\times \mathcal{H}_1(g_0)$ be the
parameter space. We define the driving system  $\theta_t$:
$\Sigma\to \Sigma$   by
\begin{equation}\label{5.1}
\theta_t(f_0(\cdot),\,g_0(\cdot))=(f_0(t+\cdot),\,g_0(t+\cdot)),
\;\; t\in \mathbb{R}.
\end{equation}

Then, system \eqref{1.1}-\eqref{1.3} is rewritten as the following
system
\begin{equation}\label{5.2}
\begin{cases}
u_{tt}+h(u_t)-\Delta u+f(u,t+s)=g(x,\,t+s), (x,t)\in \Omega\times \mathbb{R}^+,\\
u(x,\,t)|_{\partial \Omega}=0, \\
u(x,0)=u_{0}(x),\quad u_t(x,0)=v_0(x),
\end{cases}
\end{equation}
where $s\in \mathbb{R}$ means the initial symbol, corresponding to
some $\sigma \in \Sigma$.

Applying monotone operator theory or Faedo-Galerkin method, e.g.,
see \cite{CV,Li,Sh}, it is known that conditions
\eqref{1.4}-\eqref{1.9} guarantee the existence and uniqueness of
strong solution and generalized solution for
\eqref{1.1}-\eqref{1.3}, and the time-dependent terms make no
essential complications.

\begin{lemma} \label{l5.1} (\textbf{Well-posedness}) \\
Let $\Omega$ be a bounded subset of $\mathbb{R}^3$ with smooth
boundary, and assume that either Assumption  I or Assumption  II
holds. Then the non-autonomous system \eqref{5.2} has a unique
solution $(u(t),\,u_t(t))\in
\mathcal{C}(\mathbb{R}^+;\,H_0^1(\Omega)\times L^2(\Omega))$ and
$\partial_t^2u(t)\in L^2_{loc}(\mathbb{R}^+;\,H^{-1}(\Omega))$ for
any initial data $x_0=(u^{0},\,u^{1})\in H_0^1(\Omega)\times
L^2(\Omega)$ and any initial symbol $\sigma\in \Sigma$.
\end{lemma}

By $Lemma$ \ref{l5.1}, we can define the cocycle as follows:
\begin{equation}\label{5.3}
\begin{cases}
\phi:\quad \mathbb{R}^+\times \Sigma\times X & \to X, \\
(t,\sigma,(u^{0}(x),u^{1}(x))) & \to
(u^{\sigma}(t),u^{\sigma}_t(t)),
\end{cases}
\end{equation}
where $(u^{\sigma}(t),u^{\sigma}_t(t))$ is the solution of
\eqref{1.1} corresponding to initial data $(u^{0}(x),u^{1}(x))$ and
symbol $\sigma=(f_0(s+\cdot),g_0(s+\cdot))$; and for each
$(t,\,\sigma)\in \mathbb{R}^+\times \Sigma$, the mapping
$\phi(t,\,\sigma;\,\cdot): X\to X$ is continuous.

 \emph{Hereafter, we always denote  by $(\phi,\,\theta)$ the
cocycle defined in \eqref{5.1} and \eqref{5.3}.}

We now prove the following main result.

\begin{theorem}\label{t5.2} (\textbf{Existence of pullback
attractor})\\
Let $\Omega$ be a bounded domain of $\mathbb{R}^3$ with smooth
boundary. Then under either Assumption   I or Assumption  II, the
NDS $(\phi,\,\theta)$ generated by the weak solutions of
\eqref{1.1}-\eqref{1.3} has a pullback attractor
$\mathscr{A}=\{\omega_{\sigma}(B_{\sigma})\}_{\sigma \in \Sigma}$.
\end{theorem}

We need a few lemmas before proving this theorem.

\subsection{Pullback absorbing sets} \label{pullback}

In the following, we deal only with the strong solutions of
\eqref{1.1}. The generalized solution case then follows easily by
a density argument. We begin   with the following existence result
on a bounded pullback absorbing set.
\begin{lemma} \label{l5.3}(\textbf{Pullback absorbing set})\\
Let $\Omega$ be a bounded domain of $\mathbb{R}^3$ with smooth
boundary. Then under either Assumption  I or Assumption II, the NDS
$(\phi,\,\theta)$ has a bounded pullback absorbing set
$\mathscr{B}=\{B_{\sigma}\}_{\sigma\in \Sigma}$.
\end{lemma}
{\bf Proof.} For each $\sigma \in \Sigma$, we know that $\sigma$
is corresponding to some $s_0$ satisfying that
$\sigma=(f(v,s_0+t),g(x,s_0+t))$, and
$\phi(t,\theta_{-t_0}(\sigma);\,x_0)$ is the solution of the
following equation at time $t$:
\begin{equation}\label{5.4}
\begin{cases}
u_{tt}+h(u_t)-\Delta u+f(u,t-t_0+s_0)=g(x,\,t-t_0+s_0),\quad (x,t)\in \Omega\times \mathbb{R}^+,\\
(u(0),\,u_t(0))=x_0,\\
u|_{\partial \Omega}=0.
\end{cases}
\end{equation}

Under {\it Assumption $II$}, we can repeat what have done in the
proof of [$Theorem$ 1, Haraux\cite{Ha2}], to obtain that there exist
a $\rho$ (which depends only on
$\|g\|_{_{L^{\infty}(\mathbb{R},L^2(\Omega))}}$ and the coefficients
in \eqref{1.4}-\eqref{1.9}) and a $T$ (which depends only on
$\|g\|_{_{L^{\infty}(\mathbb{R},L^2(\Omega))}}$, the coefficients in
\eqref{1.4}-\eqref{1.9} and the radius of $B$) such that for any
$\sigma\in \Sigma$,
\begin{equation}\label{5.5}
\|\phi(t,\,\theta_{-t_0}(\sigma);\,x_0)\|_{_X} \leqslant \rho \quad
\text{for all}~T\leqslant t \leqslant t_0~\text{and}~x_0 \in B.
\end{equation}
Hence, for {\it Assumption $II$}, we can take $B_{\sigma}\equiv
\{x\in X~|~\|x\|_{_X} \leqslant \rho\}$ for each $\sigma$.

Under {\it Assumption $I$}, we can use the methods as that in the
proof of Chepyzhov and Vishik [\cite{CV}, $Lemma$ 4.1, $Proposition$
4.2, P 121-123], obtain also that there exist $C_{\beta}$ and
$\beta$ (which depend only on the coefficients in
\eqref{1.4}-\eqref{1.9} and $C_0$) and a $T$ (which depends only on
$\int_{-\infty}^{s_0}e^{\beta}\int_{\Omega}|g(x,s)|^2dxds$, the
coefficients in \eqref{1.4}-\eqref{1.9} and the radius of $B$) such
that
\begin{equation}\label{5.6}
\|\phi(t,\,\theta_{-t_0}(\sigma);\,x_0)\|^2_{_X} \leqslant
\rho_{_{\sigma,\beta}}=C_{\beta}(1+e^{-\beta
s_0}\int_{-\infty}^{s_0}e^{\beta s}\int_{\Omega}|g(x,s)|^2dxds)\quad
\forall~T\leqslant t \leqslant t_0,~x_0 \in B.
\end{equation}
Therefore, under {\it Assumption $I$}, we can take
$B_{\sigma}=B^{\beta}_{\sigma}= \{x\in X~|~\|x\|^2_{_X} \leqslant
\rho_{_{\sigma,\beta}}\}$. $\hfill \blacksquare$

\begin{remark}\label{r5.4}
From \eqref{5.5} we know that under {\it Assumption $II$}, the NDS
$(\phi,\,\theta)$ has a nested bounded pullback absorbing set
$\mathscr{B}=\{B_{\sigma}\}_{\sigma\in \Sigma}$.
\end{remark}

\subsection{Pullback asymptotic compactness}

We new prove  the pullback asymptotic compactness.

\begin{lemma} \label{l5.5} (\textbf{Pullback asymptotic compactness}) \\
Under either Assumption  $I$ or $II$,   for any bounded sequence
$\{x_n\}_{n=1}^{\infty}\in \mathcal{B}$ and $\sigma\in \Sigma$,
the sequence
$\phi(t_n,\,\theta_{-t_n}(\sigma);\,x_n\}_{n=1}^{\infty}$ is
precompact in $X$.
\end{lemma}
The idea for the proof is similar to that in Chueshov \& Lasiecka
\cite{CL1,CL2,CL3} and Khanmamedov \cite{Kh}; see also in
\cite{SYZ2} for linear damping and autonomous cases.


In order to prove this lemma on pullback asymptotic compactness, we
need to derive  a few energy inequalities; see
\eqref{5.19}-\eqref{5.22} below.
\medskip

We   first present some preliminaries and notations.

For each $\sigma \in \Sigma$, we know that $\sigma$ is corresponding
to some $s_0$ such that $\sigma=(f(v,s_0+t),g(x,s_0+t))$. For any
$x_0^i=(u_0^i,\,v_0^i)\in X$ ($i=1,2$), let
$(u_i(t),u_{i_t}(t))=\phi(t,\theta_{-t_0}(\sigma);\,x_0^i)$ be the
corresponding solution of the following equation at time $t$:
\begin{equation}\label{5.7}
\begin{cases}
u_{tt}+h(u_t)-\Delta u+f(u,t-t_0+s_0)=g(x,\,t-t_0+s_0),\quad (x,t)\in \Omega\times \mathbb{R}^+,\\
(u(0),\,u_t(0))=x_0^i,\\
u|_{\partial \Omega}=0.
\end{cases}
\end{equation}
For convenience, we introduce notations
\begin{equation*}
f_i(t)=f(u_i(t),\,t-t_0+s_0),\quad h_i(t)=h(u_{i_t}(t)),\quad
t\geqslant 0, ~~ i=1,2,
\end{equation*}
and
\begin{equation*}
w(t)=u_1(t)-u_2(t).
\end{equation*}

Then  $w(t)$ satisfies
\begin{equation}\label{5.8}
\begin{cases}
w_{tt}+h_1(t)-h_2(t)-\Delta w+ f_1(t)-f_2(t)=0,\\
w|_{\partial \Omega}=0,\\
(w(0),\,w_t(0))=(u_0^1,\,v_0^1)-(u_0^2,\,v_0^2).
\end{cases}
\end{equation}

We also define an energy functional
\begin{equation}\label{5.9}
E_w(t)=\frac 12\int_{\Omega}|w(t)|^2+\frac 12 \int_{\Omega}|\nabla
w(t)|^2.
\end{equation}

Since the pullback attractors obtained in $Lemma$ \ref{l5.3} are
different for \emph{Assumption I} and $II$, in the following we will
deduce different estimations.

We first deal with the case corresponding to \emph{Assumption II}:

\emph{Step 1} Multiplying \eqref{5.8} by $w_t(t)$, and integrating
over $[s,\,T]\times \Omega$, we obtain
\begin{align}\label{5.10}
E_w(T) +\int_s^T\int_{\Omega}(h_1(\tau)-h_2(\tau))w_t(\tau)dxd\tau
+\int_s^T\int_{\Omega}(f_1(\tau)-f_2(\tau))w_t(\tau)dxd\tau =
E_w(s),
\end{align}
where $0\leqslant s \leqslant T \leqslant t_0$. Then
\begin{align}\label{5.11}
\int_s^T  \int_{\Omega}(h_1(\tau)-h_2(\tau))w_t(\tau)dxd\tau
\leqslant E_w(s)-
\int_s^T\int_{\Omega}(f_1(\tau)-f_2(\tau))w_t(\tau)dxd\tau.
\end{align}
Combining with $Lemma$ \ref{l2.3}, we get that  for any $\delta
>0$,
\begin{align}\label{5.12}
\int_s^T\int_{\Omega}|w_t(\tau)|^2dxd\tau \leqslant |T-s|\delta
\cdot mes(\Omega) +C_{\delta}E_w(s) -C_{\delta}\int_s^T
\int_{\Omega}(f_1-f_2)w_txd\tau.
\end{align}

\emph{Step 2} Multiplying \eqref{5.8} by $w(t)$, and integrating
over $[0,\,T]\times \Omega$, we get that
\begin{align}\label{5.13}
\int_0^T & \int_{\Omega}|\nabla w(s)|^2dxds + \int_{\Omega}w_t(T)\cdot w(T) \nonumber \\
& = \int_0^T\int_{\Omega}|w_t(s)|^2dxds -
\int_0^T\int_{\Omega}(h_1-h_2)w + \int_{\Omega}w_t(0)\cdot w(0)  -
\int_0^T\int_{\Omega}(f_1-f_2)w.
\end{align}

Therefore, from \eqref{5.12} and \eqref{5.13}, we have
\begin{align}\label{5.14}
2\int_0^T&E_w(s)ds \nonumber \\
& \leqslant 2\delta T mes(\Omega)+2 C_{\delta}E_w(0) - 2 C_{\delta}\int_0^T\int_{\Omega}(f_1-f_2)w(t) dxds \nonumber \\
& \quad - \int_{\Omega}w_t(T)w(T) + \int_{\Omega}w_t(0)w(0)-
\int_0^T\int_{\Omega}(h_1-h_2)w-\int_0^T\int_{\Omega}(f_1-f_2)w.
\end{align}

Integrating \eqref{5.10} over $[0,\,T]$ with respect to $s$, we have
that
\begin{align}\label{5.15}
TE_w(T)&+\int_0^T\int_s^T\int_{\Omega}(h_1(\tau)-h_2(\tau))w_t(\tau) dxd\tau ds \nonumber \\
&=-\int_0^T\int_s^T\int_{\Omega}(f_1-f_2)w_tdxd\tau ds + \int_0^TE_w(s) ds \nonumber \\
& \leqslant - \int_0^T\int_s^T\int_{\Omega}(f_1-f_2)w_tdxd\tau ds+ \delta T mes(\Omega)+C_{\delta}E_w(0)  \nonumber \\
& \quad -C_{\delta}\int_0^T\int_{\Omega}(f_1-f_2)w_tdxds -\frac 12 \int_{\Omega}w_t(T)w(T) +\frac 12\int_{\Omega}w_t(0)w(0) \nonumber \\
& \quad -\frac 12 \int_0^T\int_{\Omega}(h_1-h_2)w-\frac12 \int_0^T
\int_{\Omega}(f_1-f_2)w.
\end{align}

\emph{Step 3} We will deal with $\int_0^T\int_{\Omega}(h_1-h_2)w$.
Multiplying \eqref{5.7} by $u_{i_t}(t)$,   we obtain
\begin{align*}
\frac 12\frac{d}{dt}\int_{\Omega}(|u_{i_t}|^2+|\nabla
u_i|^2)+\int_{\Omega}h(u_{i_t})u_{i_t}+\int_{\Omega}f(u_i,t+s_i)u_{i_t}=\int_{\Omega}g_iu_{i_t},
\end{align*}
which, combining with the existence of bounded uniformly absorbing
set, implies that
\begin{align}\label{5.16}
\int_0^T\int_{\Omega}h(u_{i_t})u_{i_t} \leqslant M_T,
\end{align}
where the constant $M_T$ depends on $T$ (which is different from the
autonomous cases). Then, noticing \eqref{2.1}, we obtain that
\begin{equation}\label{5.17}
\int_0^T\int_{\Omega}|h(u_{i_t})|^{\frac{p+1}{p}}dxds \leqslant M_T.
\end{equation}

Therefore, using H$\ddot{o}$lder inequality, from \eqref{5.17} we
have
\begin{align*}
|\int_0^T\int_{\Omega}h_iw|  \leqslant
M_T^{\frac{p}{p+1}}\left(\int_0^T\int_{\Omega}|w|^{p+1}
\right)^{\frac{1}{p+1}},
\end{align*}
which implies that
\begin{align}\label{5.18}
|\int_0^T\int_{\Omega}(h_1-h_2)w|  \leqslant
2M_T^{\frac{p}{p+1}}\left(\int_0^T\int_{\Omega}|w|^{p+1}
\right)^{\frac{1}{p+1}}.
\end{align}

Hence, combining \eqref{5.15} and \eqref{5.18}, we obtain that
\begin{align*}
E_w(T)& \leqslant \delta  mes(\Omega)-\frac{1}{T} \int_0^T\int_s^T\int_{\Omega}(f_1(\tau)-f_2(\tau))w_t(\tau)dxd\tau ds+ \frac{C_{\delta}}{T}E_w(0)  \nonumber \\
& \quad -\frac{C_{\delta}}{T}\int_0^T\int_{\Omega}(f_1(s)-f_2(s))w_t(s)dxds -\frac{1}{2T} \int_{\Omega}w_t(T)w(T) +\frac{1}{2T}\int_{\Omega}w_t(0)w(0) \nonumber \\
& \quad +\frac 1T
M_T^{\frac{p}{p+1}}\left(\int_0^T\int_{\Omega}|w(s)|^{p+1}dxds
\right)^{\frac{1}{p+1}}-\frac{1}{2T} \int_0^T
\int_{\Omega}(f_1(s)-f_2(s))w(s)dxds
\end{align*}
for any $0\leqslant T\leqslant t_0$.

We define
\begin{align}\label{5.19}
&\psi_{T,\,\delta,\,\sigma}(x_0^1,\,x_0^2) \nonumber \\
&=-\frac{1}{T} \int_0^T\int_s^T\int_{\Omega}(f_1(\tau)-f_2(\tau))w_t(\tau)dxd\tau ds  \nonumber \\
& \quad -\frac{C_{\delta}}{T}\int_0^T\int_{\Omega}(f_1(s)-f_2(s))w_t(s)dxds -\frac{1}{2T} \int_{\Omega}w_t(T)w(T)  \nonumber \\
& \quad +\frac 1T
M_T^{\frac{p}{p+1}}\left(\int_0^T\int_{\Omega}|w(s)|^{p+1}dxds
\right)^{\frac{1}{p+1}}-\frac{1}{2T} \int_0^T
\int_{\Omega}(f_1(s)-f_2(s))w(s)dxds.
\end{align}
Then we   have
\begin{equation}\label{5.20}
E_w(T) \leqslant \delta mes(\Omega)
+\frac{1}{2T}\int_{\Omega}w_t(0)w(0)+
\frac{C_{\delta}}{T}E_w(0)+\psi_{T,\,\delta,\,\sigma}(x_0^1,\,x_0^2)
\end{equation}
for any $\delta >0$, $0\leqslant T\leqslant t_0$.

\medskip

For the case corresponding to \emph{Assumption I}:

Since under our general assumption \eqref{1.10}, as shown in
\eqref{5.6}, the pullback attractors may not satisfy the nested
properties. Inspired partly by the results in \cite{CLR,LLR} we will
deduce different estimations by the same methods; see \eqref{5.21}
and \eqref{5.22} below.

Repeat \emph{Step 1} and \emph{Step 2} above, and just replace the
multipliers $w_t(t)$ and $w(t)$ by $e^{\beta t}w_t(t)$ and $e^{\beta
t}w(t)$ respectively, and take into account $\beta<C_0$, we can
obtain the following similar estimates
\begin{equation}\label{5.21}
E_w(T) \leqslant \frac{\alpha}{TC_0}e^{-\beta
T}E_w(0)+\psi'_{T,\,\sigma}(x_0^1,\,x_0^2)
\end{equation}
for any $\delta >0$, $0\leqslant T\leqslant t_0$, where
$\alpha=(C_0+\beta)/(C_0-\beta)$ and
\begin{align}\label{5.22}
&\psi'_{T,\,\sigma}(x_0^1,\,x_0^2)
\nonumber \\
&=-\frac{e^{-T\beta}}{T} \int_0^T\int_s^T\int_{\Omega}e^{\beta
\tau}(f_1(\tau)-f_2(\tau))w_t(\tau)dxd\tau ds
+\frac{\alpha}{2T}e^{-\beta
T}\int_{\Omega}w_t(0)w(0)\nonumber \\
& \quad -\frac{\alpha}{TC_0}e^{-\beta
T}\int_0^T\int_{\Omega}e^{\beta s}(f_1(s)-f_2(s))w_t(s)dxds -\frac{C}{T} \int_{\Omega}w_t(T)w(T)  \nonumber \\
& \quad + M_{_{E_w(0), T,
C_0,\beta}}\left(\int_0^T\int_{\Omega}|w(s)|^{2}dxds
\right)^{\frac{1}{2}}.
\end{align}

With the above energy inequalities, we are now
ready to prove pullback asymptotic compactness. \\

{\bf Proof of $Lemma$ \ref{l5.5}:} We will deal with
\emph{Assumption I} and \emph{Assumption II} separately.

\emph{Assumption II}:

For each $\sigma\in \Sigma$, and for any fixed $\varepsilon>0$, from
\eqref{5.19}, we can take $t_0$ large enough such that
\begin{equation}
E_w(t_0) \leqslant
\varepsilon+\psi_{t_0,\,\delta,\,\sigma}(x_0^1,\,x_0^2)~\text{for
all}~x_0^1,\,x_0^2\in B_{\sigma}.
\end{equation}
Hence, thanks to $Theorem$ \ref{t4.2} and $Lemma$ \ref{l5.3}, it is
sufficiently to prove that the function
$\psi_{t_0,\,\delta,\,\sigma}(\cdot,\,\cdot)$ defined in
\eqref{5.19} belongs to $Contr(B_{\sigma})$ for each fixed $t_0$.

We observe from equation \eqref{5.7} (and also see \cite{Ha1}) that
for any $t_0>0$,
\begin{equation}\label{5.24}
\bigcup_{t \in
[0,\,t_0]}\phi(t,\,\theta_{-t_0}(\sigma);B_{\sigma})~~\text{is
bounded in}~ X,
\end{equation}
and the bound depends only on $t_0$ and $\sigma$.

Let $(u_n,u_{t_n})$ be the corresponding solution of
$(u^n_0,v_0^n)\in B_{\sigma}$ for problem \eqref{5.7},
$n=1,2,\cdots$. From the observation above, without loss of
generality (or by passing to subsequences), we assume that
\begin{equation}\label{5.25}
u_n \to u\quad  \star-\text{weakly in} ~ L^{\infty}(0,t_0;\,
H_0^1(\Omega)),
\end{equation}
\begin{equation}\label{5.26}
u_n \to u\quad  \text{in} ~L^{p+1}(0,t_0;\, L^{p+1}(\Omega)),
\end{equation}
\begin{equation}\label{5.27}
u_{n_t} \to u_t\quad  \star-\text{weakly in} ~ L^{\infty}(0,t_0;\,
L^2(\Omega)),
\end{equation}
\begin{equation}\label{5.28}
u_n \to u\quad  \text{in} ~ L^2(0,t_0;\, L^2(\Omega))
\end{equation}
and
\begin{equation}\label{5.29}
u_n(0) \to u(0)~~\text{and}~~u_n(t_0) \to u(t_0)\quad  \text{in} ~
L^4(\Omega).
\end{equation}
Here we have used the compact embeddings $H_0^1 \hookrightarrow
L^{4}$ and $H_0^1 \hookrightarrow L^{p+1}$ (since $1 \leq p<5$).

Now, we will deal with each term in \eqref{5.19} one by one.

First, from \eqref{5.24}, \eqref{5.29} and \eqref{5.26} we get that
\begin{equation}\label{5.30}
\lim_{n\to \infty}\lim_{m\to \infty}\int_{\Omega}
(u_{n_t}(t_0)-u_{m_t}(t_0))(u_n(t_0)-u_m(t_0)) dx = 0,
\end{equation}
\begin{equation}\label{5.31}
\lim_{n\to \infty}\lim_{m\to
\infty}\int_0^{t_0}\int_{\Omega}|u_n(s)-u_m(s)|^{p+1}dxds=0,
\end{equation}
and from \eqref{1.8} and \eqref{5.28}, we further have
\begin{equation}\label{5.32}
\lim_{n\to \infty}\lim_{m\to \infty}\int_0^{t_0}\int_{\Omega}
(f(u_n(s),s-t_0+s_0)-f(u_m(s),s-t_0+s_0))(u_n(s)-u_m(s))dxds = 0.
\end{equation}

Second, note that
\begin{align*}
\int_0^{t_0}\int_{\Omega}&
(u_{n_t}(s)-u_{m_t}(s))(f(u_n(s),s-t_0+s_0)-f(u_m(s),s-t_0+s_0))dxds \nonumber \\
&  = \int_0^{t_0}\int_{\Omega}
u_{n_t}(s)f(u_n(s),s-t_0+s_0)+\int_0^{t_0}\int_{\Omega}
u_{m_t}(s)f(u_m(s),s-t_0+s_0)\nonumber \\
&\qquad -\int_0^{t_0}\int_{\Omega}  u_{n_t}(s)f(u_m(s),s-t_0+s_0)-\int_0^{t_0}\int_{\Omega} u_{m_t}(s)f(u_n(s),s-t_0+s_0)\nonumber \\
& =\int_{\Omega} F(u_n(t_0),s_0)-\int_{\Omega}
F(u_n(0),-t_0+s_0)-\int_0^{t_0}\int_{\Omega}F_s(u_n(\tau),\tau-t_0+s_0)dxd\tau
\nonumber \\
&\qquad+\int_{\Omega} F(u_m(t_0),s_0)-\int_{\Omega}
F(u_m(0),-t_0+s_0) -\int_0^{t_0}\int_{\Omega}F_s(u_m(\tau),\tau-t_0+s_0)dxd\tau \nonumber \\
&\qquad -\int_0^{t_0}\int_{\Omega}
u_{n_t}(s)f(u_m(s),s-t_0+s_0)-\int_0^{t_0}\int_{\Omega}
u_{m_t}(s)f(u_n(s)s-t_0+s_0).
\end{align*}
By   \eqref{5.25}, \eqref{5.27}, \eqref{5.29} and \eqref{1.8},
taking first $m\to \infty$ and then $n\to \infty$, we obtain that
\begin{align}\label{5.33}
\lim_{n\to \infty}&\lim_{m\to \infty}\int_0^{t_0}\int_{\Omega}
(u_{n_t}(s)-u_{m_t}(s))(f(u_n(s),s-t_0+s_0)-f(u_m(s),s-t_0+s_0))dxds \nonumber \\
&=\int_{\Omega} F(u(t_0),s_0)-\int_{\Omega}
F(u(0),-t_0+s_0)-\int_0^{t_0}\int_{\Omega}F_s(u(\tau),\tau-t_0+s_0)dxd\tau
\nonumber \\
&\qquad+\int_{\Omega} F(u(t_0),s_0)-\int_{\Omega}
F(u(0),-t_0+s_0) -\int_0^{t_0}\int_{\Omega}F_s(u(\tau),\tau-t_0+s_0)dxd\tau \nonumber \\
&\qquad -\int_0^{t_0}\int_{\Omega}
u_tf(u(s),s-t_0+s_0)-\int_0^{t_0}\int_{\Omega} u_{t}f(u(s),s-t_0+s_0) \nonumber \\
&=0.
\end{align}
Similarly, we have
\begin{align*}
\int_s^{t_0}\int_{\Omega}&
(u_{n_t}(\tau)-u_{m_t}(\tau))(f(u_n(\tau),\tau-t_0+s_0)-f(u_m(\tau),\tau-t_0+s_0))dxd\tau \nonumber \\
& =\int_{\Omega} F(u_n(t_0),s_0)-\int_{\Omega}
F(u_n(s),s-t_0+s_0)-\int_s^{t_0}\int_{\Omega}F_s(u_n(\tau),\tau-t_0+s_0)dxd\tau
\nonumber \\
&\quad+\int_{\Omega} F(u_m(t_0),s_0)-\int_{\Omega}
F(u_m(s),s-t_0+s_0) -\int_s^{t_0}\int_{\Omega}F_s(u_m(\tau),\tau-t_0+s_0)dxd\tau \nonumber \\
&\quad -\int_s^{t_0}\int_{\Omega}
u_{n_t}f(u_m(\tau),\tau-t_0+s_0)-\int_s^{t_0}\int_{\Omega}
u_{m_t}f(u_n(\tau),\tau-t_0+s_0).
\end{align*}
Since $|\int_s^{t_0}\int_{\Omega}
(u_{n_t}(\tau)-u_{m_t}(\tau))(f(u_n(\tau),\tau-t_0+s_0)-f(u_m(\tau),\tau-t_0+s_0))dxd\tau|$
is bounded for each fixed $t_0$,   by the Lebesgue dominated
convergence theorem, we finally have
\begin{align}\label{5.34}
\lim_{n\to \infty}&\lim_{m\to
\infty}\int_0^{t_0}\int_s^{t_0}\int_{\Omega}
(u_{n_t}(\tau)-u_{m_t}(\tau))(f(u_n(\tau),\tau-t_0+s_0)-f(u_m(\tau),\tau-t_0+s_0))dxd\tau ds \nonumber \\
&= \int_0^{t_0}\left(\lim_{m\to \infty}\lim_{n\to
\infty}\int_s^{t_0}\int_{\Omega}
(u_{n_t}(\tau)-u_{n_t}(\tau))\right.
\nonumber \\
&\qquad \qquad \qquad \qquad \qquad \qquad
(f(u_n(\tau),\tau-t_0+s_0)-f(u_m(\tau),\tau-t_0+s_0))dxd\tau
\Big)ds \nonumber \\
&=\int_0^{t_0} 0ds=0.
\end{align}

Hence, from \eqref{5.30}-\eqref{5.34}, we see that
$\psi_{t_0,\delta,\sigma}(\cdot,\,\cdot)\in Contr(B_{\sigma})$.

\emph{Assumption I}:

From the definition of $\rho_{_{\sigma,\beta}}$ (see \eqref{5.6}),
we have the conclusion: {\it for each $\sigma$ and for any
$\varepsilon
>0$, we can take $t_0$ large enough such that $e^{-\beta
t_0}\rho_{_{\theta_{-t_0}(\sigma),\beta}}\leqslant \varepsilon$}.

Hence, from $Corollary$ \ref{c4.3} and \eqref{5.21}, we only need to
verify that the function $\psi'_{t_0,\,\sigma}(\cdot,\,\cdot)$
defined in \eqref{5.22} belongs to
$Contr(B_{\theta_{-t_0}(\sigma)})$. To this end, we notice that
$e^{\beta t}$ is bounded in $[0,t_0]$ and $\bigcup_{t \in
[0,\,t_0]}\phi(t,\,\theta_{-t_0}(\sigma);B_{\theta_{-t_0}(\sigma)})$
is bounded in $X$. The remainder is just a repeat of that for
$\psi_{t_0,\,\delta,\,\sigma}(\cdot,\,\cdot)$ above.

This completes the proof of $Lemma$ \ref{l5.5}.
 $\hfill \blacksquare$

\subsection{Existence of pullback attractors}

Now we complete the proof of the main result.

 {\bf Proof of
$Theorem$ \ref{t5.2}} \ From $Lemma$ \ref{l5.3} and $Lemma$
\ref{l5.5}, we see that the conditions of $Theorems$ \ref{t3.13} and
\ref{t3.14} are all satisfied respectively and thus we imply the
existence of the pullback attractor. $\hfill \blacksquare$

\begin{remark} \label{kappa}
In this section, we obtain the pullback asymptotic compactness for
the non-autonomous wave system \eqref{1.1}-\eqref{1.3} by the
technique  presented in \S 4.  This technique is different from the
method in \cite{CLR}. Due to the existence of nested bounded
pullback absorbing set for \emph{Assumption II}(Lemma \ref{l5.3}),
the pullback $\kappa-$contraction is equivalent to pullback
asymptotic compactness (see $Theorem$ \ref{t3.13}). Thus, in
principle, we could also use the pullback $\kappa-$contraction
criterion to conclude the existence of pullback attractor, using the
decomposition method as in \cite{Fe,SYZ1} (popular for autonomous
systems).

\end{remark}

\section{Some remarks}

In this paper, we discuss the asymptotic behavior of solutions in
the framework of pullback attractors. Another interesting question
is forward attractors; see \cite{Ch,CKS} for general discussions
or \cite{CKR} for practical applications to wave equations with
delays. However, for the forward attraction property to hold, one
usually needs some uniformity  about the time-dependent terms
(i.e., about   the symbol spaces \cite{CV}).  As discussed in
details in \cite{CLR1,LLR},   for   general non-autonomous
dissipative systems,  how to obtain the forward attraction
properties is an open problem if without this uniformity
assumption.

For our problem, under the \textbf{Assumption  II}   in \S 1, we
indeed obtain a bounded uniformly absorbing set in the sense of
\cite{CV} in $Lemma$ \ref{l5.3} (or see Haraux\cite{Ha1}). If we
assume further that $g$ satisfies some additional conditions, e.g.,
$g$ is translation compact or $g\in
W^{1,\infty}(\mathbb{R};\,L^2(\Omega))$, then by the same method, we
can verify the family of processes (see \cite{CV} for more details)
corresponding to the non-autonomous wave system
\eqref{1.1}-\eqref{1.3} is uniformly asymptotically compact and thus
has a uniform (w.r.t. $\sigma \in \Sigma$) attractor in the sense of
\cite{CV}. However, for the case of \textbf{Assumption I}
  in \S 1, it appears difficult to discuss the forward
attraction  for $g$ satisfying only \eqref{1.10}.

For the \emph{autonomous} case  of \eqref{1.1}-\eqref{1.3},
recently, Chueshov \& Lasiecka \cite{CL3} have shown a   general
result   for the existence of global attractor, and they allow $p
=5$, i.e., the so-called critical interior damping. In their
\emph{autonomous} case, it is true that for all $0\leqslant
s\leqslant t$,
\begin{equation}\label{6.1}
\int_s^t\int_{\Omega}h(u_t)u_tdxd\tau \leqslant C_R,
\end{equation}
where $C_R$  depends only on the norm of initial data, but
independent of time instants $s$ and $t$. However, for our
non-autonomous case, this constant may depend on time instants $s$
and $t$ (e.g., see \eqref{5.16},\eqref{5.17}), and thus in our
proofs, we require (at least, technically) that the growth order of
$h$ to be strictly less than $5$: $p<5$.

Moreover, in the present paper, we use the {\it pullback asymptotic
compactness} to obtain the existence of pullback attractors of
non-autonomous hyperbolic systems. This is mainly based on a
technical method for verifying {\it pullback asymptotic compactness}
in \S 4. However, for other non-autonomous systems
or using other techniques (e.g., the   decomposition method), the
pullback $\kappa-$contraction criterion may be more appropriate
for proving the existence of pullback attractors.

Finally, we point out that all the contents in this paper can be
expressed by the framework of \emph{processes}, instead of
cocycles, as   in  \cite{CKR,CR,CLR1,LLR}; see also \cite{CV}
  for more results about processes.
\\

\noindent{\bf \Large Acknowledgement}

The authors would like to thank Edriss Titi and  the referees  for
helpful comments and suggestions.


\begin{thebibliography}{99}

\bibitem{ACH} J. Arrieta, A.  N. Carvalho and J. K. Hale, \textit{A damped
hyperbolic equations with critical exponents}, Comm. Partial
Differential Equations, 17(1992) 841-866.

\vspace{-0.3cm}

\bibitem{ARS} B. Aulbach, M. Rasmussen and S. Siegmund, \textit{Approximation of
attractors of nonautonomous dynamical systems}, Discrete Contin.
Dyn. Syst. Ser. B, 5(2005) 215-238.

\vspace{-0.3cm}

\bibitem{BV} A. V. Babin and M. I. Vishik, \textit{Attractors of evolution equations}, North-Holland, Amsterdam,
1992.

\vspace{-0.3cm}

\bibitem{Ba} J. M. Ball, \textit{Global attractors for damped semilinear wave equations},
Discrete Contin. Dyn. Syst., 10(2004) 31-52.

\vspace{-0.3cm}

\bibitem{BCD} V. P. Bongolan-Walsh, D. Cheban and J. Duan, \textit{Recurrent motions in
the nonautonomous Navier-Stokes system}, Discrete Contin. Dyn. Syst.
B, 3(2002) 255-262.

\vspace{-0.3cm}

\bibitem{CKR} T. Caraballo, P. E. Kloeden and J. Real,
\textit{Pullback and forward attractors for a damped wave equation
with delays}, Stoch. Dyn., 4(2004) 405-423.

\vspace{-0.3cm}

\bibitem{CLR} T. Caraballo, G. Lukaszewicz and J. Real,
\textit{Pullback attractors for asymptotically compact
non-autonomous dynamical systems}, Nonlinear Analysis, TMA, 64(2006)
484-498.

\vspace{-0.3cm}

\bibitem{CLR1} T. Caraballo, G. Lukaszewicz and J.Real, \textit{Pullback
attractors for non-autonomous 2D-Navier-Stokes equations in some
unbounded domains}, C.R. Acad. Sci. Paris, Ser.I 342(2006)
263-268.

\vspace{-0.3cm}

\bibitem{CR} T. Caraballo and J. Real, \textit{Attractors for
2D-Navier-Stokes models with delays}, J. Differential Equations,
205(2004) 271-297.

\vspace{-0.3cm}

\bibitem{Cha} S. Chandrasekhar, \textit{The mathematical theory of
black holes}, Oxford University Press, 1992.

\vspace{-0.3cm}

\bibitem{Ch} D. N. Cheban, \textit{Global attractors of non-autonomous dissipative dynamical systems},
World Scientific Publishing, 2004.

\vspace{-0.3cm}

\bibitem{CD} D. N. Cheban and J. Duan, \textit{Almost periodic motions and global attractors
of the non-autonomous Navier-Stokes equations}, J. Dyn. Differential
Equations, 16(2004) 1-34.

\vspace{-0.3cm}

\bibitem{CKS} D. N. Cheban, P.E. Kloeden and B. Schmalfuss, \textit{The relationship between pullback, forwards
and global attractors of nonautonomous dynamical systems}, Nonlinear
Dyn. Syst. Theory, 2(2002) 9-28.

\vspace{-0.3cm}

\bibitem{CV} V. V. Chepyzhov and M.I. Vishik, \textit{Attractors for equations of mathematical physics}, Amer. Math. Soc. Colloq. Publ., Vol. 49, Amer. Math. Soc.,
Providence, RI, 2002.

\vspace{-0.3cm}

\bibitem{CY} L. Chierchia and J. You, \textit{KAM tori for 1D wave equations
with periodic boundary conditions}, Comm. Math. Phys., 211(2000)
497-525.

\vspace{-0.3cm}

\bibitem{CL1} I. Chueshov and I. Lasiecka, \textit{Attractors for second-order evolution equations
with nonlinear damping}, J. Dyn. Differential Equations, 16(2004)
469-512.

\vspace{-0.3cm}

\bibitem{CL2} I. Chueshov and I. Lasiecka, \textit{Long-time behavior of second evolution
equations with nonlinear damping}, Memoirs AMS, to appear, 2006.

\vspace{-0.3cm}

\bibitem{CL3} I. Chueshov and I. Lasiecka, \textit{Long time
dynamics of semilinear wave equation with nonlinear
interior-boundary damping and sources of critical exponents}, to
appear.

\vspace{-0.3cm}

\bibitem{Fe} E. Feireisl, \textit{Global attractors for damped wave equations with
supercritical exponent}, J. Differential Equations, 116(1995)
431-447.

\vspace{-0.3cm}

\bibitem{GM} J. M. Ghidagla and A. Marzocchi, \textit{Longtime behavior of strongly damped wave equations,
global attractors and their dimension}, SIAM J. Math. Anal.,
22(1991) 879-895.

\vspace{-0.3cm}

\bibitem{Ha} J. K. Hale, \textit{Asymptotic behavior of dissipative systems}, Amer. Math. Soc., Providence, RJ, 1988.

\vspace{-0.3cm}

\bibitem{Ha1} A. Haraux, \textit{Two remarks on dissipative hyperbolic problems}, dans "Nonlinear
partial differential equations and their applications, College de
France Seminar", vol. 7 (H. Brezis \& J.L. Lions editors), Research
Notes in Math., Pitman, 122(1984) 161-179.

\vspace{-0.3cm}

\bibitem{Ha2} A. Haraux, \textit{Recent results on semilinear wave equations
with dissipation}, Pitman Research Notes in Math., 141(1986)
150-157.

\vspace{-0.3cm}

\bibitem{Ilyin} A. A. Ilyin, A. Miranville and E. S. Titi, \textit{Small viscosity sharp
estimates for the global attractor of the 2-D damped-driven
Navier-Stokes equations}, Commun. Math. Sci., 2(2004) 403-426.

\vspace{-0.3cm}

\bibitem{Kh} A. Kh. Khanmamedov, \textit{Global
attractors for von Karman equations with nonlinear interior
dissipation}, J. Math. Anal. Appl., 318(2006) 92-101.

\vspace{-0.3cm}

\bibitem{LLR} J. Langa, G. Lukaszewicz and J.Real, \textit{Finite fractal
dimension of pullback attractors for non-autonomous 2D Navier-Stokes
equations in some unbounded domains}, Nonlinear Analysis, TMA,
dio:10.1016/j.na.2005.12.017.

\vspace{-0.3cm}

\bibitem{LSC} P. S. Landahl, O. H. Soerensen and P. L. Christiansen,
\textit{Soliton excitations in Josephson tunnel junctions}, Phys.
Rev. B, 25(1982) 5737-5748.

\vspace{-0.3cm}

\bibitem{LS1} J. A. Langa and B. Schmalfuss, \textit{Finite dimensionality of
attractors for non-autonomous dynamical systems given by partial
differential equations}, Stoch. Dyn., 3(2004) 385-404.

\vspace{-0.3cm}

\bibitem{Li} J. L. Lions, \textit{Quelques methodes de resolution des
problemes aux limites nlineaires}, Dunod, Paris, 1969.

\vspace{-0.3cm}

\bibitem{LS} G. Lukaszewicz and W. Sadowski, \textit{Uniform attractor for 2D
magneto-micropolar fluid flow in some unbounded domains}, Z. Angew.
Math. Phys., 55(2004) 1-11.

\vspace{-0.3cm}

\bibitem{Majda} A. Majda, \textit{Introduction to PDEs and waves for the
atmosphere and oceans}, American Math. Soc., Providence, RI, 2003.

\vspace{-0.3cm}

\bibitem{MRW} I. Moise, R. Rosa and X. Wang, \textit{Attractors for noncompact
nonautonomous systems via energy equations}, Discrete Contin. Dyn.
Syst., 10(2004) 473-496.

\vspace{-0.3cm}

\bibitem{PZ} V. Pata and S. Zelik, \textit{A remark on the damped wave
equation}, Comm. Pure Appl. Anal., 5(2006) 609-614.

\vspace{-0.3cm}

\bibitem{Pedlosky} J. Pedlosky, \textit{Geophysical Fluid Dynamics},
Springer-Verlag, 2nd edition, 1987.

\vspace{-0.3cm}

\bibitem{Robinson} J.C. Robinson, \textit{Infinite-dimensional dynamical
systems}, Cambridge University Press, 2001.

\vspace{-0.3cm}

\bibitem{SY} G. R. Sell and Y. You, \textit{Dynamics of evolutionary equations},
Springer, New York, 2002.

\vspace{-0.3cm}

\bibitem{Sh} R. Showalter, \textit{Monotone operators in Banach
spaces and nonlinear partial differential equations}, Amer. Math.
Soc., Providence, RI, 1997.

\vspace{-0.3cm}

\bibitem{SYZ1} C. Y. Sun, M. H. Yang and C. K. Zhong, \textit{Global attractors for
the wave equation with nonlinear damping}, J. Differential
Equations, 227(2006) 427-443.

\vspace{-0.3cm}

\bibitem{SYZ2} C. Y. Sun, M. H. Yang and C. K. Zhong, \textit{Global attractors for hyperbolic equations with critical exponent in
locally uniform spaces}, submitted.

\vspace{-0.3cm}

\bibitem{Te} R. Temam, \textit{Infinite-dimensional dynamical systems in mechanics and physics}, Springer-Verlag New York, 1997.

\vspace{-0.3cm}

\bibitem{ZYS} C. K. Zhong, M. H. Yang and C. Y. Sun, \textit{The existence of global attractors for the
norm-to-weak continuous semigroup and its application to the
nonlinear reaction-diffusion equations}, J. Differential Equations,
223(2006), 367-399.

\vspace{-0.3cm}

\bibitem{ZW} S. F. Zhou and L. S. Wang, \textit{Kernel sections for damped
non-autonomous wave equations with critical exponent}, Discrete
Contin. Dyn. Syst., 9(2003) 399-412.
\end{thebibliography}
\end{document}